\documentclass[reqno]{amsart}

\usepackage[mathscr]{eucal}

\usepackage{amsfonts, amsmath, amsthm, amssymb, latexsym, amsbsy}
\newtheorem{theorem}{Theorem}[section]
\newtheorem{lemma}[theorem]{Lemma}
\newtheorem{prop}[theorem]{Proposition}
\newtheorem{cor}[theorem]{Corollary}
\theoremstyle{definition}

\newtheorem{example}[theorem]{Example}
\theoremstyle{remark}
\newtheorem{remark}[theorem]{Remark}

\numberwithin{equation}{section}

\newcommand{\GL}{\mathtt{GL}}

\newcommand{\Spec}{\mathrm{Spec}}

\newcommand{\ad}{\mathrm{ad}}

\newcommand{\End}{\mathrm{End}}
\newcommand{\Hom}{\mathrm{Hom}}

\newcommand{\id}{1}

\def\C{\mathbb{C}}

\def\CM{\mathcal{C}}

\def\Z{\mathbb{Z}}

\def\Irr{\mathtt{Irr}}

\def\N{\mathbb{N}}
\def\P{\mathcal{P}}
\def\wx{\widehat X}
\def\wy{\widehat Y}
\def\llambda{{\boldsymbol{\lambda}}}
\def\aalpha{{\boldsymbol{\alpha}}}

\def\h{\mathfrak{h}}

\def\gr{\mbox{\rm{Gr}}^{\mbox{\scriptsize{\rm{ad}}}}}
\def\sl2{{\mathfrak{s}\mathfrak{l}}_2}

\def\tr{\mathrm{tr}}
\newcommand{\ee}{{\mathfrak{e}}}
\newcommand{\diag}{\mathop{\mathrm{diag}}\nolimits}

\usepackage{hyperref}

\usepackage[utf8]{inputenc}
\usepackage[english]{babel}

\def\Rep{{\mathtt{Rep}}}

\newcommand{\p}{\partial}

\numberwithin{equation}{section}

\renewcommand{\le}{\leqslant}
\renewcommand{\ge}{\geqslant}
\renewcommand{\cong}{\simeq}

\newcommand{\wt}{\widetilde}
\newcommand{\wh}{\widehat}

\newcommand{\ra}{\rangle}

\newcommand{\HH}{\hbox to 15pt {${\cal H}$\hss${\cal H}$}}

\newcommand{\hreg}{\mathfrak h_{\mathrm reg}}

\begin{document}

\title{KP hierarchy for the cyclic quiver}

\author{Oleg Chalykh}
\address{School of Mathematics, University of Leeds, Leeds LS2 9JT, UK}
\email{o.chalykh@leeds.ac.uk}
\author{Alexey Silantyev}
\address{School of Mathematics, University of Leeds, Leeds LS2 9JT, UK}
\email{aleksejsilantjev@gmail.com}
%

%
\maketitle
\begin{abstract}
We introduce a generalisation of the KP hierarchy, closely related to the cyclic quiver and the Cherednik algebra $H_k(\Z_m)$. This hierarchy depends on $m$ parameters (one of which can be eliminated), with the usual KP hierarchy corresponding to the $m=1$ case. Generalising the result of G. Wilson~\cite{W1}, we show that our hierarchy admits solutions parameterised by suitable quiver varieties. The pole dynamics for these solutions is shown to be governed by the classical Calogero--Moser system for the wreath-product $\Z_m\wr S_n$ and its new spin version. These results are further extended to the case of the multi-component hierarchy.
\end{abstract}

\section{Introduction}
\label{intro}

It is a remarkable fact, first observed in~\cite{AMM}, that integrable non-linear PDEs admit solutions whose poles move as classical-mechanical particles in a completely integrable fashion. One of the best-known examples is the KP hierarchy, in which case the classical integrable system in question is the Calogero--Moser system, see~\cite{K, CC}. Since the Calogero--Moser system can be generalised to all root systems and Coxeter groups~\cite{OP1,OP2}, there is a natural question whether these generalisations can be similarly linked to integrable non-linear equations. For the infinite series of classical root systems a partial answer can be obtained by folding~\cite{OT,BST}, which relates the Calogero--Moser systems in types $A$ and $B$. However, what can be obtained in type $B$ by folding constitutes only a very special case, while for the most general coupling parameters the answer, to the best of our knowledge, remained unknown. The present paper gives an answer to this question for all classical root systems and much more: as we will explain, there is a natural generalisation of the KP hierarchy that admits rational solutions whose pole dynamics is governed by the Calogero--Moser system for the complex reflection group $\Z_m\wr S_n$. The case $m=1$ is the usual, type $A$ Calogero--Moser problem, and the case $m=2$ covers the remaining classical types.

Our approach to this problem is inspired by a result of G. Wilson, who has shown in~\cite{W1} that rational solutions of the KP hierarchy can be parametrised by (the points of) certain affine varieties, the Calogero--Moser spaces. The $n$-th Calogero--Moser space $\CM_n$ can be viewed as a completed phase space of the $n$-particle Calogero--Moser system (hence their name). On the other hand, $\CM_n$ is a particular example of a quiver variety, with the underlying quiver consisting of just one loop. From that perspective, the present work is a generalisation of that of Wilson, with the one-loop quiver replaced by the cyclic quiver on $m$ vertices. It is worth mentioning that the spaces $\CM_n$ have also appeared in the context of non-commutative projective geometry and classification of ideals of the first Weyl algebra~\cite{BW1, BW2, KKO}. The results of~\cite{BW1} prompted Crawley-Boevey and Holland to conjecture a similar link in the case of affine Dynkin quivers (see~\cite[p.~45]{CB1}). Their conjecture was later generalised and proved by Baranovsky, Ginzburg and Kuznetsov in~\cite{BGK1, BGK2}. In the case of the cyclic quiver this is related to the problem of classifying ideals of the Cherednik spherical subalgebra  $\textbf{e}H_k\textbf{e}$ for the cyclic group. In~\cite{Esh} this problem was tackled differently to~\cite{BGK1, BGK2}, by using the approach of~\cite{BC}. In fact, some of our formulas and calculations closely resemble those in~\cite{BC} and~\cite{Esh}. 

While the moduli spaces parametrising solutions to our hierarchy are known from the above works, the hierarchy itself seems new. We construct it in a Lax form in the same spirit as the usual KP hierarchy, but replacing the operator $\p=d/dx$ with the one-variable Dunkl operator for $\Gamma=\Z_m$. Solutions to the hierarchy for given quiver data are described by an explicit formula, so checking that they satisfy the equations of the hierarchy amounts to a direct algebraic calculation. Note that we do not require a passage to an adelic Grassmannian $\gr$ as in~\cite{W1} (for the details of how to define the appropriate notion in the present context see~\cite{BHY, BGK2}). It is worth mentioning that our initial motivation for the present work was to understand the moduli space of bispectral rings of differential operators constructed by Bakalov, Horozov and Yakimov in~\cite{BHY}. An attempt in that direction, also inspired by~\cite{W1}, was made by Rothstein in~\cite{Ro}, but the proofs in {\it loc. cit.} were omitted. Moreover, as our work shows, the answer to that question turns out to be simpler and more natural than it might be suggested by~\cite{Ro}. Apart from a brief remark in the concluding section, however, we do not discuss this problem in detail here, leaving it to a separate publication.

Our generalisation of the KP hierarchy is naturally related not just to a cyclic quiver but to a suitable framing of it, obtained by adding a vertex and a number of arrows. Different choices of the framing lead to different hierarchies and to different integrable multi-particle systems governing the pole dynamics. In the simplest case of one added arrow the classical-mechanical system in question can be identified with the rational Calogero--Moser system with the symmetry group $G=\Z_m\wr S_n$, while the other cases naturally lead to spin versions of that system. In the case $m=1$ one recovers the Gibbons--Hermsen system \cite{GH}, while in general the systems we obtain this way seem new. Note that there exist various versions of the spin Calogero--Moser system associated to simple complex Lie algebras, see \cite{LX} and references therein. However, to the best of our knowledge, spin versions of the Calogero--Moser system for $G=\Z_m\wr S_n$ and $m>2$ did not appear in the literature previously. Because we define these systems as a commutative family of Hamiltonian flows on suitable Calogero--Moser spaces, the fact that the flows are complete is obvious from the start. As another by-product we also find a Lax matrix for the (spin) Calogero--Moser system for $G=\Z_m\wr S_n$. Note that $G$ is not a Coxeter group for $m>2$, so the corresponding Calogero--Moser system and its Lax matrix can not be obtained within the general framework of \cite{BCS}.

The plan of the paper is as follows. In Section \ref{sec2} we recall the definition of the Cherednik algebra for $\Gamma=\Z_m$, use its microlocalisation to construct our version of the KP hierarchy in a Lax form, and introduce its spherical sub-hierarchy. In Section \ref{sec3} we recall the definition of deformed preprojective algebras for an arbitrary quiver $Q$ due to~\cite{CBH} and define quiver varieties associated to a framing of $Q$. We spell this out in detail for a cyclic quiver $Q_0$. By analogy with \cite{W1}, one calls these varieties Calogero--Moser spaces. In Section \ref{secRatSol} we obtain our main results: we introduce a natural Hamiltonian dynamical system on a Calogero--Moser space and construct rational solutions of our hierarchy starting from any point of this space. A similar result is obtained in Section \ref{secSolSubh} for the spherical sub-hierarchy. The integrable systems governing the pole dynamics for the constructed solutions are considered in more detail in Section~\ref{sec5}. For that we first recall a link, due to Etingof and Ginzburg~\cite{EG}, between the Calogero--Moser spaces and Cherednik algebras for $G=\Z_m\wr S_n$. We use this link to identify the flows of the spherical sub-hierarchy with the Calogero--Moser system for the wreath-product $\Z_m\wr S_n$. The particle system related to the full hierarchy is considered in Section \ref{secFullHier}; it can be viewed as a spin version of the Calogero--Moser system for $\Z_m\wr S_n$, where the spin variables belong to the group algebra $\C\Z_m$. Section~\ref{secMultGen} describes the multi-component version of our hierarchy; this extends the results of the previous two sections and leads to a construction of the general spin version of the Calogero--Moser system for $G=\Z_m\wr S_n$. The concluding section discusses some open questions, further directions and links.

{\it Acknowledgement}. We are grateful to Yu.~Berest, W.~Crawley-Boevey, P.~Etingof and P.~Iliev for stimulating discussions, and to G.~Bellamy for useful comments. Both authors gratefully acknowledge the support by EPSRC under grant EP/K004999/1.

\section{Cherednik algebra for the cyclic group and generalised KP hierarchy}
\label{sec2}

\subsection{The Cherednik algebra of $\Gamma=\Z_m$.}

Let $\Gamma=\{e, {s}, \dots, {s}^{m-1}\}$ be the multiplicative cyclic group of order $m$ with generator ${s}$. Identifying $\Gamma$ with the group $\Z_m$ of roots of unity, we have an action of $\Gamma$ on $\C$, with ${s}$ acting as the multiplication by a fixed primitive $m$th root of unity, $\mu$. This induces an action of $\Gamma$ on the ring $\C(x)$ of rational functions and differential operators of one complex variable. Explicitly, ${s}(f)(x)=f(\mu^{-1}x)$ for $f\in\C(x)$, and ${s}(\p)=\mu\p$, $\p:=\frac{d}{d x}$. Therefore, we can form the crossed product $\mathcal D(x, \p)*\Gamma$, where $\mathcal D(x, \p)$ denotes the ring of differential operators in $x$ with rational coefficients. Recall that for a group $G$ acting on an algebra $A$, the crossed product $A*G$ is spanned by the pairs $(a,g)$ with $a\in A$, $g\in G$, with the multiplication defined by $(a,g)(b,h)=(ag(b), gh)$. Abbreviating $(a,g)$ as $ag$, the multiplication consists of moving the group elements to the right using the rule $ga=g(a)g$ so that $(ag)(bh)=ag(b)gh$. For instance, in our situation we have
\begin{align}\label{cr}
 {s} x =\mu^{-1} x {s}\qquad\text{ and}\qquad {s}\p=\mu\p {s}\,.
\end{align} 
The differential filtration on the ring $\mathcal D(x, \p)$ extends to a filtration on $\mathcal D(x, \p)*\Gamma$, by placing the elements of $\Gamma$ in degree $0$.

The group algebra $\C\Gamma$ has a natural basis given by the idempotents
\begin{align}\label{id}
 &\epsilon_i=\frac1m\sum_{r=0}^{m-1}\mu^{-ir}{s}^r\in\C[\Gamma] &&(i=0,\dots, m-1)\,.
\end{align}
Given $k=(k_0, \dots, k_{m-1})\in\C^m$ and $\tau\ne 0$, one defines the Cherednik algebra $H_{\tau,k}(\Gamma)$ as the subalgebra of $\mathcal D(x, \p)*\Gamma$, generated by ${s}$, $x$ and the following element, called {\it Dunkl operator}:
\begin{align}\label{du}
y=\tau\p-\frac{1}{x}\sum_{i=0}^{m-1}mk_i\epsilon_i\,. 
\end{align}
It is clear from~\eqref{cr} that ${s} y=\mu y{s}$ and that
\begin{align}
 \epsilon_ix=x\epsilon_{i+1}\,,\qquad \epsilon_iy=y\epsilon_{i-1}\,. \label{epsilonxy}
\end{align}
(Here and in similar situations below the indices are always taken modulo $m$.)
Also, the following relation between $x$ and $y$ is easy to check:
\begin{align}\label{cr1} 
 &xy-yx=\sum_{i=0}^{m-1}\lambda_i\epsilon_i\,,&&\text{with}\quad \lambda_i=m(k_{i-1}-k_{i})-{\tau}\,.
\end{align}
It follows that the algebra $H_{\tau,k}$ can be described by generators and relations as
\begin{align}\label{is}
 &H_{\tau,k}\cong (\C\langle x, y\ra*\Gamma)/(xy-yx-\lambda)\,,&&\lambda=\sum_{i=0}^{m-1}\lambda_i\epsilon_i\in\C\Gamma\,,
\end{align} 
where $\C\langle x, y\ra$ denotes the free algebra on $x, y$ and the action of $\Gamma$ is given by ${s}(x)=\mu^{-1}x$, ${s}(y)=\mu y$.

The algebra $H_{\tau,k}$ is a particular example of a rational Cherednik algebra, which can be associated to any complex reflection group (see~\cite{EG}). It has a PBW property, namely, a vector space isomorphism $H_{\tau,k}\cong \C[x]\otimes\C\Gamma\otimes\C[y]$, which means that every element can be uniquely presented as a linear combination of the terms $x^i\gamma y^j$ with $i,j\ge 0$ and $\gamma\in\Gamma$. The inclusion $H_{\tau,k}\subset \mathcal D(x, \p)*\Gamma$ induces a filtration on $H_{\tau,k}$, called {\it differential filtration}, with $\deg y=1$, $\deg x=0$, $\deg\gamma=0$.

The algebra~\eqref{is} has also appeared in the context of non-commutative deformations of Kleinian singularities~\cite{Ho, CBH}, and it is isomorphic to the deformed preprojective algebra of a cyclic quiver as defined in~\cite{CBH} (we will recall that definition and isomorphism in Section~\ref{sec3} below). Note that if we rescale $\tau$ and $k_i$ simultaneously, the resulting Cherednik algebra will not change. 
For $\tau=0$ one defines $H_{0,k}$ by replacing the Dunkl embedding with its quasi-classical counterpart, see~\cite{EG}; this agrees with the isomorphism~\eqref{is}.

\subsection{Microlocalisation of $H_{\tau,k}$} 
The well-known construction of the KP hierarchy uses the ring of pseudo-differential operators $\mathcal D(x,\p)((\p^{-1}))$. Our generalisation is based instead on a microlocalisation $(\mathcal D(x, \p)*\Gamma)((y^{-1}))$, where $y$ is the Dunkl operator~\eqref{du}. To construct such a microlocalisation explicitly, we will first consider the Ore localisation of $H_{\tau,k}$ at $y$ and then pass to a completion.

We begin by introducing the parameters $c_\gamma$, $\gamma\in\Gamma$ determined from 
\begin{align}\label{cgamma}
\sum_{i=0}^{m-1}mk_i\epsilon_i=\sum_{\gamma\in\Gamma}c_\gamma\gamma\,.
\end{align}
Explicitly, $c_\gamma=\sum_{i=0}^{m-1}\mu^{-ij}k_i$ for $\gamma={{s}^j}$. The following relation in $H_{\tau,k}$ is immediate from~\eqref{du} and~\eqref{cgamma}:
\begin{align*}
ya-ay=\tau\p(a)-\sum_{\gamma\in\Gamma\setminus\{e\}}c_\gamma x^{-1}(\gamma a -a\gamma)\qquad\text{for any}\ a=a(x)\in\C[x]\,. 
\end{align*}
A similar commutation relation holds for any $f\in\C[x]*\Gamma$. Namely, let us define an algebra automorphism $\theta$ of $\C[x]*\Gamma$ by $\theta({s})=\mu{s}$, $\theta(x)=x$. In other words,
\begin{align*}  
\theta\bigg(\sum_{i=0}^{m-1} a_i(x){s}^i\bigg)=\sum_{i=0}^{m-1} a_i(x)\mu^i{s}^i\,.
\end{align*}
Clearly, $\theta^m=\id$. With this notation, one has
\begin{align}\label{D}
 &yf-\theta^{-1}(f)y=\tau\p(f)-\sum_{\gamma\in\Gamma\setminus\{e\}}c_\gamma x^{-1}(\gamma f -f\gamma) &&\text{for any}\ f\in\C[x]*\Gamma\,. 
\end{align}
Here $\p$ acts on $\C[x]*\Gamma$ by the formula $\p\big(\sum_i a_i(x){s}^i\big)=\sum_i \p a_i(x){s}^i$. Note that for $a=\sum_i a_i(x){s}^i$ we have
\begin{align*}
 x^{-1}({s}^j a- a{s}^j)=\sum_i x^{-1}\big(a_i(\mu^{-j}x)-a_i(x)\big){s}^{i+j}\,, 
\end{align*}
therefore, the denominators in~\eqref{D} cancel, as expected. The right-hand side of~\eqref{D} defines a linear map $D\colon\C[x]*\Gamma\to\C[x]*\Gamma$, that is,
\begin{align}\label{DD}
D(f)=\tau\p(f)-\sum_{\gamma\in\Gamma\setminus\{e\}}c_\gamma x^{-1}[\gamma, f]\,.
\end{align}
From~\eqref{D} we have $yf=\theta^{-1}(f)y+D(f)$. By induction, for any $n\in\N$,  
\begin{align}\label{form}
y^n f=\sum_{p=0}^n \sum_{\genfrac{}{}{0pt}{1}{i_0, \dots, i_{p}\ge 0}{i_0+\dots+i_p=n-p}}\theta^{-i_0}D\theta^{-i_1}\dots D\theta^{-i_p}(f)y^{n-p}\,.
\end{align}
Note that an application of $D$ to a polynomial in $\C[x]*\Gamma$ lowers its degree, therefore the right-hand side of~\eqref{form} will terminate at $p=\deg(f)$. This implies that the set $\{y^k\,|\,k\in\N\}$ is (right and left) Ore, so let $H_{\tau,k}[y^{-1}]$ be the corresponding right localisation. By the PBW property of $H_{\tau,k}$, every element of $H_{\tau,k}[y^{-1}]$ has a unique presentation as a finite sum
$u=\sum_{i\in\Z} u_i y^i$ with $u_i\in\C[x]*\Gamma$. The explicit form of the multiplication in $H_{\tau,k}[y^{-1}]$ can be derived from~\eqref{form}. Namely, one first finds that
\begin{align}\label{rule}
y^{-1}f=\sum_{p=0}^\infty (-1)^p(\theta D)^p\theta(f)y^{-1-p}\,.
\end{align}
From that, using induction we get
\begin{align}\label{rule1}
y^{-n} f=\sum_{p=0}^\infty (-1)^p\sum_{\genfrac{}{}{0pt}{1}{\varepsilon_1+\dots+\varepsilon_{n+p-1}=p}{\varepsilon_i\in \{0, 1\}}}\theta D^{\varepsilon_1}\theta D^{\varepsilon_2}\dots D^{\varepsilon_{n+p-1}}\theta(f)y^{-n-p}\,.
\end{align}

Note that if $f$ is not polynomial but, say, rational, $f\in \C(x)*\Gamma$, then the sum in~\eqref{rule} and~\eqref{rule1} is no longer finite. Therefore, in analogy with the notion of formal pseudo-differential operators, let us introduce
\begin{align}
\label{p}
\mathcal P=\left\{\sum_{i\ge 0}f_j y^{N-j}\ \mid\ f_j\in\C(x)*\Gamma\,,\ N\in\Z\right\}\,. 
\end{align}  
Addition in $\mathcal P$ is obvious, while the multiplication is determined by the rules~\eqref{form}--\eqref{rule1}. The fact that the multiplication is associative follows from its associativity in $H_{\tau, k}[y^{-1}]$. Indeed, checking the associativity in every degree in $y$ requires checking an algebraic relation involving a finite number of the coefficients $f_j\in\C(x)*\Gamma$ and operations $\theta$ and $D$. Since we know that this relations hold for arbitrary polynomial coefficients $f_j$, they must hold identically. Hence, the multiplication in $\mathcal P$ is associative; note that it depends on the parameters $\tau, k$ entering into the definition of $D$.

Fixing $N$ in the definition of $\mathcal P$, one gets a subspace $\mathcal P_N\subset\mathcal P$ and an increasing chain $\dots\subset \mathcal P_N\subset\mathcal P_{N+1}\subset\dots$, which makes $\mathcal P$ into a filtered ring. For any 
$S=\sum_{j\ge 0}f_jy^{N-j}$ in  $\mathcal P$, we have $S=S_-+S_-$ where
\begin{align}
 &S_-=\sum_{j > N}f_j\,y^{N-j}\,, &&S_+= \sum_{j=0}^N f_j\,y^{N-j}\,. \label{SmSp}
\end{align}
(In the case $N<0$ we put $S_-=S$ and $S_+=0$.) This defines a direct sum decomposition $\mathcal P=\mathcal P_-\oplus\mathcal P_+$ with 
\begin{align*}
 &\mathcal P_-=\mathcal P_{-1}=\Bigg\{\sum_{j \ge 1} f_j y^{-j} \Bigg\} &&\text{and} &&\mathcal P_+=\Bigg\{\sum_{j=0}^N f_j y^{j}\ \mid\ \ N\ge 0\Bigg\}\,.
\end{align*}
From the formula for $y$ it easily follows that $\mathcal P_+=\mathcal D(x, \p)*\Gamma$. The following properties of the multiplication in the ring $\mathcal P$ should be clear:
\begin{equation*}
(Sf)_{\pm}=S_{\pm}\cdot f\qquad\forall S\in\mathcal P\,,\ \forall f\in\C(x)*\Gamma\,.
\end{equation*}

\begin{remark}\label{rem1} An alternative construction of the ring $\mathcal P$ is by applying microlocalisation in the sense of~\cite{AVdBVO}. Indeed, $\mathcal P$ is isomorphic to the microlocalisation $\mathcal Q^\mu_S$ of $R=\mathcal D(x, \p)*\Gamma$ at the multiplicatively closed set $S=\{y^k\,|\,k\in\N_0\}$. This can be checked either at the level of the associated graded rings, using Propositions~3.5 and~3.10 of {\it loc. cit.}, or by direct comparison.  
\end{remark}

\begin{remark}\label{rem2}
In the definition of $\mathcal P$ one can set the coefficients $f_j$ to be from any ring $\mathcal F*\Gamma$, closed under the action of $D$. For instance, we may replace $\mathcal F=\C(x)$ by $\C[[x]]$, $\C((x))$ or by the rings of functions, analytic (respectively, meromorphic) near $x=0$. We may further replace any such $\mathcal F$ by a ring of matrices $\mathrm{Mat}(n, \mathcal F)$.   
\end{remark}

\begin{remark}\label{rem3}
The formulas~\eqref{form},~\eqref{rule1} can be rewritten using the operations $D_i:=\theta^{-i+1}D\theta^{i}$. Then for any $n\in\N$ one has:
\begin{gather}\label{form2}
y^n\theta^n(f)=\sum_{\varepsilon_1, \dots, \varepsilon_n\in\{0,1\}}D_1^{\varepsilon_1}\dots D_n^{\varepsilon_n}(f)y^{n-\varepsilon_1-\dots -\varepsilon_n}\,,\\
\label{rule2}
y^{-n}\theta^{-n}(f)=\sum_{p=0}^\infty(-1)^p\sum_{\genfrac{}{}{0pt}{1}{l_1,\dots, l_n\ge 0}{l_1+ \dots  +l_n=p}} (D_0)^{l_n}(D_{-1})^{l_{n-1}}\dots (D_{1-n})^{l_1}(f)\, y^{-n-p}\,.
\end{gather}

\end{remark}

\subsection{Generalised KP hierarchy}

Let $\mathcal P$ be as in Remark~\ref{rem2} above, i.e. we allow arbitrary coefficents $f_j\in\mathcal F*\Gamma$, where $\mathcal F$ is either the ring $\C(x)$ of rational functions or one of its variants. 
Consider an element $L\in \mathcal P$ of the form
\begin{align}
\label{L}
L=y+\sum_{i\ge 0}f_i y^{-i}\,,
\end{align}
where the coefficients $f_i\in\mathcal F*\Gamma$ depend on an infinite number of additional variables $t_1, t_2, \dots$, called times.
Our generalisation of the KP hierarchy is the following system of equations:
\begin{align}
 &\frac{\partial L}{\partial t_k}=\big[(L^{k})_+, L\big]\,,&&k=1,2,\dots\,.
\label{kp} 
\end{align}
Here $(L^k)_+$ denotes the positive part of $L^k$ according to the notation~\eqref{SmSp}, 
 and $[\cdot , \cdot]$ is the commutator bracket, $[A, B]=AB-BA$. The system~\eqref{kp} is viewed as a hierarchy of non-linear equations on the infinite number of dependent variables $f_i(x)=f_i(x; t_1, t_2, \dots)\in\mathcal F*\Gamma$. When $m=1$ we have $\Gamma=\{e\}$ and $y=\partial$, in which case the system~\eqref{kp} turns into the standard presentation of the KP hierarchy.

Since $L\in y+\mathcal P_0$, each equation in~\eqref{kp} can be viewed as a flow on $\mathcal P_0$. The following results are analogous to the similar standard facts about the KP hierarchy.
\begin{prop} \label{prop2.1}
{\normalfont (1)} The equations~\eqref{kp} imply the zero-curvature equations:
\begin{align}
 &\Big[\frac{\partial}{\partial t_k}-(L^k)_+\;,\;\frac{\partial}{\partial t_l}-(L^l)_+\Big]=0 &&\forall\ k,l\in \N\,. \label{zcc}
\end{align}

{\normalfont (2)} The system~\eqref{kp} is consistent, that is, the equations of the hierarchy define commuting flows on $\mathcal P$.

\end{prop}

Proofs of both statements are the same as for the usual KP hierarchy, see, e.g., \cite[Sec. 1.6]{Di}. \qed

\medskip

For the usual KP hierarchy, it is customary to assume that $f_0=0$ in~\eqref{L}. The analogous condition in our case is as follows. For any $f\in\mathcal F*\Gamma$, denote by $f_\gamma\in\mathcal F$ the coefficient in the decomposition $f=\sum_{\gamma\in\Gamma}f_\gamma \gamma$; in particular, $f_e$ denotes the coefficient at the identity element. Then our assumption on the coefficient $f_0$ in~\eqref{L} will be as follows:
\begin{align}
\label{f0}
(f_0)_e=0\,.
\end{align} 
To see that this is preserved by the flows of the hierarchy, rewrite the equations~\eqref{kp} as
\begin{align}
\label{kpminus}
 \frac{\partial L}{\partial t_k}=-\big[(L^{k})_-,L\big]\,.
\end{align}
Writing $(L^{k})_-$ as $hy^{-1}+\dots$, we obtain that $\frac{\partial f_0}{\partial t_k}= \theta^{-1}(h)-h$, from which $\frac{\partial (f_0)_e}{\partial t_k}=(\theta^{-1}(h)-h)_e=0$, as needed.

Let us recall that for the KP hierarchy ($m=1$), assuming $f_0=0$, one has $L_{+}=\partial$ and so the first equation in \eqref{kp} gives simply that $\frac{\partial L}{\partial t_1}=\partial(L)$, which means that one can set $t_1=x$. Since in our case $f_0$ is nonzero, the first equation becomes more complicated and we can no longer identify $x$ with one of the time variables. Still, the zero-curvature equations \eqref{zcc} for all $k, l$ less than some given $n$ lead to a closed set of nonlinear equations on a finite number of coefficients of $L$. For example, to get an analogue of the usual KP equation, one can take the zero-curvature equations for all $k, l\in\{1,2,3\}$. However, the calculations become too cumbersome to be carried out explicitly.    

Our main goal will be to construct some explicit solutions to this hierarchy. The solutions will be constructed in the form
\begin{align}
\label{dress}
 &L=MyM^{-1}\,, &&\text{where}\quad  M\in 1+\mathcal P_-\,.
\end{align}
Writing $M$ as $M=1+\sum_{i\ge1}w_i(x)y^{-i}$, we can calculate $M^{-1}$, term by term, from the equation $M^{-1}M=1$: 
\begin{align*}
M^{-1}=1-w_1y^{-1}+(w_1\theta(w_1)-w_2)y^{-2}+\dots \,.
\end{align*} 
Substituting these into~\eqref{dress} gives
\begin{align*}
L=y-(\theta^{-1}(w_1)-w_1)+(-D(w_1)+w_2-\theta^{-1}(w_2)-w_1^2+\theta^{-1}(w_1)w_1)y^{-1}+\dots \,.
\end{align*}
We therefore can view~\eqref{dress} as a non-linear transformation between the variables $f_j$ and $w_i$, so that 
\begin{align*}
f_0=-\theta^{-1}(w_1)+w_1\,,\quad f_1=-D(w_1)+w_2-\theta^{-1}(w_2)-w_1^2+\theta^{-1}(w_1)w_1\,,
\end{align*}
and so on. It is clear from the first formula that~\eqref{f0} is automatically satisfied.

\begin{remark}
The above $M^{-1}$ is the left inverse to $M$. In the same manner one can construct the right inverse. By a standard argument, they coincide so $M^{-1}$ is the two-sided inverse.  
\end{remark}

The following result is straightforward.
\begin{prop} \label{wave}
Suppose $M\in 1+\mathcal P_-$ satisfies the equation
\begin{align*}
 \frac{\partial M}{\partial t_k}=-(L^{k})_- M
\end{align*}
for some $k$. Then $L=MyM^{-1}$ satisfies the equation
\begin{align*}
\frac{\partial L}{\partial t_k}=-\big[(L^{k})_-, L\big]\,.
\end{align*}  
\end{prop}

\subsection{Spherical sub-hierarchy}
\label{spherical}

Now suppose that $L$ is as above and in addition it satisfies the following property with respect to the action of the group $\Gamma$:
\begin{align}
\label{sym} 
\epsilon_i L= L \epsilon_{i-1}\,.
\end{align}
Then it is easy to see that such condition is only compatible with the equations~\eqref{kp} when $k$ is divisible by $m$. Therefore, in this case we have a sub-hierarchy consisting of the equations
\begin{align}
 &\frac{\partial L}{\partial t_{mj}}=\big[(L^{mj})_+, L\big]\,, &&j\in\N\,.
\label{kpm} 
\end{align}  
It is convenient to rewrite this as equations on $\widetilde L=L^m$. It follows from~\eqref{kpm} that 
\begin{align}
\label{kpmm}
 &\frac{\partial \widetilde L}{\partial t_{mj}}=\big[(\widetilde L^{j})_+, \widetilde L\big]\,, &&j\in\N\,.
\end{align}
We have
\begin{align}
\label{Ltilde}
 &\widetilde L=y^m+\sum_{a > 0}g_ay^{m-1-a}\,, &&\epsilon_i\widetilde L= \widetilde L\epsilon_i\,, &&g_a\in\mathcal F *\Gamma\,.
\end{align}
(The fact that $g_0=0$ is an easy corollary of~\eqref{f0}.)

We will call the hierarchy~\eqref{kpmm}--\eqref{Ltilde} the {\it spherical sub-hierarchy}. It can be further decoupled into $m$ independent hierarchies. Namely, 
decompose $\widetilde L$ as
\begin{align*}
 &\widetilde L=\sum_{\ell=0}^{m-1} \widetilde L_\ell\,, &&\widetilde L_\ell=\epsilon_\ell\widetilde L=\epsilon_\ell\widetilde L\epsilon_\ell\,.
\end{align*}
Then the system~\eqref{kpmm} splits into $m$ independent systems for each of $\widetilde L_\ell$:
\begin{align}
\label{kpmmm}
 &\frac{\partial \widetilde L_\ell}{\partial t_{mj}}=\big[(\widetilde L_\ell^{j})_+, \widetilde L_\ell\big]\,, &&\widetilde L_\ell=\epsilon_\ell y^m\epsilon_\ell+\sum_{a>0}\epsilon_\ell g_ay^{m-1-a}\epsilon_\ell\,.
\end{align}
In fact, the systems \eqref{kpmmm} with different $\ell$ are all equivalent to each other, up to a change of parameters, and if $\widetilde L_0, \dots, \widetilde L_{m-1}$ are solutions to each of these hierarchies then
$\widetilde L:=\widetilde L_0+\dots+\widetilde L_{m-1}$ is a solution to \eqref{kpmm}. 

\begin{remark} Although the cyclic quiver did not appear explicitly above, the Cherednik algebra $H_{\tau, k}(\Z_m)$ can be in fact interpreted as a deformed preprojective algebra associated to the cyclic quiver with $m$ vertices, see Section \ref{casecyclic} below. The role of quivers will become more apparent when we turn to constructing rational solutions of the hierarchy in Section \ref{secRatSol}.
\end{remark}

\section{Deformed preprojective algebras and Calogero--Moser spaces} 
\label{sec3}

We begin by defining the deformed preprojective algebras due to Crawley-Boevey and Holland, and by recalling some facts about their representations; the interested reader should consult~\cite{CBH, CB2} for the details. Consider a quiver $Q$, an oriented graph with the set of vertices $I$ and arrows $Q$. For an arrow $a\colon i\to j$, we use the notation $t(a), h(a)$ to denote the tail $i$ and the head $j$ of $a$. Consider the double quiver $\overline Q$, obtained by adjoining a reverse arrow $a^*\colon j\to i$ for every arrow $a\colon i\to j$ in $Q$. 
Let $\mathbb C\overline Q$ be the path algebra of $\overline Q$; it has basis the paths in $\overline Q$, including a trivial path $e_i$ for each vertex. We use the convention that for $a\colon i\to j$ and $b\colon j\to k$, the product $ba$ represents the path $i\to j\to k$ and that $ba=0$ if $h(a)\ne t(b)$.  

Given $\lambda\in \mathbb C^{I}$,  the {\it deformed preprojective algebra} $\Pi^\lambda(Q)$ is defined as the quotient of $\C\overline Q$ by the following relation: 
\[
\sum_{a\in Q} (aa^*-a^*a)= \sum_{i\in I} \lambda_ie_i\,.
\]
Multiplying this relation by $e_i$, we obtain an equivalent set of relations, one for each vertex:
\begin{align}
\label{rel}
 &\sum_{a\in Q\atop h(a)=i}aa^*-\sum_{a\in Q\atop t(a)=i}a^*a=\lambda_ie_i\,, &&i\in I\,. 
\end{align} 
For simplicity, we will simply write $\Pi^\lambda$ for $\Pi^\lambda(Q)$.

Representations of $\overline Q$ of dimension $\alpha\in\N_{0}^I$ are given by attaching vector spaces $V_i=\C^{\alpha_i}$ to each vertex $i\in I$, with the arrows of the quiver represented by linear maps $X_a\colon V_{t(a)} \to V_{h(a)}$ (with each trivial path $e_i$ represented by the identity map on $V_i$). The set of all representations of a given dimension is, therefore,
\begin{align*}
\mathtt{Rep}(\C\overline Q, \alpha)=\bigoplus_{a\in Q} \mathtt{Mat}(\alpha_{h(a)} \times \alpha_{t(a)}, \C) \oplus \mathtt{Mat}(\alpha_{t(a)} \times \alpha_{h(a)}, \C)\,.
\end{align*}
Isomorphism classes of representations correspond to orbits under the group 
\begin{align}\label{galpha}
G(\alpha)=\bigg(\prod_{i\in I} \mathtt{GL}(\alpha_i, \mathbb C)\bigg)/\,\C^\times\,,
\end{align}
acting by conjugation. Here $\C^\times$ denotes the diagonal subgroup of scalar matrices, acting trivially on $\mathtt{Rep}(\C\overline Q, \alpha)$. Semi-simple representations correspond to closed orbits.

We have naturally $\mathtt{Rep}(\Pi^\lambda, \alpha)\subset \mathtt{Rep}(\C\overline Q, \alpha)$ as a subvariety cut out by the relations~\eqref{rel}:
\begin{align}
\label{relx}
 &\sum_{a\in Q\atop h(a)=i}X_aX_{a^*}-\sum_{a\in Q\atop t(a)=i}X_{a^*}X_a=\lambda_i\id_{V_i}\,, &&i\in I\,. 
\end{align} 
By adding these and taking trace, we obtain 
\begin{align}
\label{tra}
\lambda\cdot\alpha=\sum_{i\in I}\lambda_i\alpha_i=0\,.
\end{align} 
Consider the quotient $\mathtt{Rep}(\Pi^\lambda, \alpha)//G(\alpha)$, whose points correspond to semi-simple representations of $\Pi^{\lambda}$ of dimension $\alpha$. We will be mostly interested in the situation where $\lambda$ and $\alpha$ are such that a general representations in $\mathtt{Rep}(\Pi^\lambda, \alpha)$ is simple. In this case, the following result is known by~\cite[Corollary~1.4 and Lemma~6.5]{CB2}.
\begin{prop} 
\label{dim}
Let $p(\alpha)=1+\sum_{a\in Q}\alpha_{t(a)}\alpha_{h(a)}-\alpha\cdot\alpha$, where $\alpha\cdot\alpha=\sum_{i\in I}\alpha_i^2$. If a general representation in $\mathtt{Rep}(\Pi^\lambda, \alpha)$ is simple, then $\mathtt{Rep}(\Pi^\lambda, \alpha)//G(\alpha)$ is a reduced and irreducible scheme of dimension $2p(\alpha)$, smooth at the points corresponding to (the isomorphism classes of) simple representations.
In particular, if all representations in $\mathtt{Rep}(\Pi^\lambda, \alpha)$ are simple, then $\mathtt{Rep}(\Pi^\lambda, \alpha)//G(\alpha)$ is a smooth connected affine variety of dimension $2p(\alpha)$. 
\end{prop}
 
Let us remark that the varieties $\mathtt{Rep}(\Pi^\lambda, \alpha)//G(\alpha)$ can be obtained by Hamiltonian reduction~\cite{MW}. Indeed, the space $\mathtt{Rep}(\C\overline Q, \alpha)$ has a natural complex symplectic structure
\begin{align}
\label{symp}
\omega=\sum_{a\in Q}\tr(\mathrm{d}X_{a^*} \wedge \mathrm{d}X_{a})\,,
\end{align}
invariant under $G(\alpha)$. The left-hand side of~\eqref{relx} can then be viewed as a component $\mu_i$ of a moment map $\mu$ for the action of $G(\alpha)$ on $\mathtt{Rep}(\C\overline Q, \alpha)$. Therefore, we can view
$\mathtt{Rep}(\Pi^{\lambda}, \alpha)$ as a fibre $\mu^{-1}(\overline\lambda)$ of the moment map at $\overline\lambda=\sum_{i\in I}\lambda_i\id_{V_i}$:
\begin{align}
\label{sympq}
\mathtt{Rep}(\Pi^\lambda, \alpha)//G(\alpha)=\mu^{-1}(\overline\lambda)//G(\alpha)\,.
\end{align}
As a consequence, in the situation when $G(\alpha)$ acts freely, the quotient~\eqref{sympq} carries a natural symplectic structure.

\subsection{Quiver varieties}
  
Let $Q_0$ be an arbitrary quiver with the vertex set $I$. A framing of $Q_0$ is a quiver $Q$ that has one additional vertex, denoted $\infty$, and a number of arrows $i\to \infty$ from the vertices of $Q_0$ (multiple arrows are allowed). Given $\alpha\in \N_0^I$ and $\lambda\in\C^I$, we extend them from $Q_0$ to $Q$ by putting $\alpha_\infty=1$ and $\lambda_\infty=-\lambda\cdot\alpha$. Thus, we put
\begin{align}
\label{frame}
\aalpha=(1, \alpha)\,, \qquad \llambda=(-\lambda\cdot\alpha, \lambda)\,.
\end{align}
This choice ensures that $\llambda\cdot\aalpha=0$. Consider the space $\mathtt{Rep}(\Pi^{{\llambda}}, \aalpha)$ of representations of the deformed preprojective algebra of $Q$ of dimension $\aalpha$. The quotients
\begin{align}
\label{qvar}
\mathfrak M_{\alpha,\lambda} (Q):=\mathtt{Rep}(\Pi^{\llambda}, \aalpha)//G(\aalpha)
\end{align} 
are closely related to Nakajima quiver varieties, see~\cite{CB2}, pp. 261--262. Note that since $\alpha_\infty=1$, we have 
\begin{align}
\label{gl}
G(\aalpha)=\bigg(\prod_{i\in I\sqcup\{\infty\}} \mathtt{GL}(\alpha_i, \mathbb C)\bigg)\big/\C^\times \cong \mathtt{GL}(\alpha):=\prod_{i\in I} \mathtt{GL}(\alpha_i, \mathbb C)\,,
\end{align}
where the group on the right-hand side acts on $\bigoplus_{i\in I} V_i$.

Recall that Kac defines in~\cite{Kac} the root system $\Delta(Q)$ for an arbitrary quiver~$Q$. Let us say that $\lambda\in\C^I$ is regular if $\lambda\cdot\alpha\ne 0$ for any $\alpha\in\Delta(Q_0)$. Then the following result is a corollary of~\cite[Theorem~1.2]{CB2}, (cf.~\cite[Proposition~3]{BCE}):

\begin{prop}
Let  $\aalpha$ and $\llambda$ be as in~\eqref{frame}, with $\lambda$ regular. Choose a framing $Q$ of $Q_0$ and let $\Pi^{\llambda}=\Pi^{\llambda}(Q)$. Then there is a $\Pi^{{\llambda}}$-module of dimension $\aalpha$ if and only if $\aalpha\in\Delta_+(Q)$, where $\Delta_+(Q)$ is the set of positive roots of $Q$. Moreover, every such module is simple. 
\end{prop}

Combining this with Prop.~\ref{dim}, we conclude that for $(1,\alpha)\in\Delta_+(Q)$ and regular $\lambda\in\C^I$ the variety $\mathfrak M_{\alpha,\lambda} (Q)$ is smooth, of dimension $2p(\aalpha)$. Note that for imaginary roots $\aalpha$ one has $p(\aalpha)>0$, while for real roots $p(\aalpha)=0$. Furthermore, if $\aalpha$ is a real positive root, then there is a unique representation of $\Pi^\llambda$ of dimension $\aalpha$ up to isomorphism~\cite[Proposition~1.2]{CB3}. This means that in this case $\mathfrak M_{\alpha,\lambda} (Q)$ consists of one point.

\subsection{Case of a cyclic quiver}\label{casecyclic}

Now let $Q_0$ be a cyclic oriented quiver with the set of vertices $I=\mathbb Z/m\mathbb Z=\{0,1, \dots, m-1\}$ and $m$ arrows $a_i\colon i+1\to i$. The double quiver $\overline Q_0$ is obtained by adding the reverse arrows $a_i^*\colon i\to i+1$. For $\lambda\in \C^I$, the deformed preprojective algebra $\Pi^\lambda$ is the quotient of $\C\overline Q_0$ by the relations
\begin{align}
\label{cycrel}
 &a_{i}a_{i}^*-a_{i-1}^*a_{i-1}=\lambda_ie_i\,, &&i\in I\,.
\end{align}
If $m=1$, then we have a one-loop quiver $Q_0$, so $\Pi^\lambda$ in this case is isomorphic to $\mathbb C\langle a, a^*\rangle/\{aa^*-a^*a=\lambda\}$. For $m>1$ there is an isomorphism between $\Pi^\lambda$ and Cherednik algebra $H_{\tau, k}(\mathbb Z_m)$, given by
\begin{align}
\label{isoH}
 a_i\mapsto \epsilon_i x \,,\qquad a_i^*\mapsto y\epsilon_i\,,\qquad e_i\mapsto \epsilon_i\,,
\end{align}
with the parameters $\lambda_i$ and $\tau, k$ related by~\eqref{cr1}.

\subsection{Calogero--Moser spaces}

Let $Q_0$ be a cyclic oriented quiver with the set of vertices $I=\mathbb Z/m\mathbb Z$ as above. Calogero--Moser spaces are the varieties~\eqref{qvar} associated to certain framings of $Q_0$. We will consider one of the following framings:
\begin{align}
\label{qd}
Q&=Q_0\sqcup\{b_{i,\alpha}\colon i\to\infty\,|\, i\in I\,,\ \alpha=1, \dots, d\}\,,  &&d\in\N\,,\\
\label{qld}
Q&=Q_0\sqcup\{b_{\ell,\alpha}\colon\ell\to\infty\,|\, \alpha=1, \dots, d\}\,, &&\ell\in I\,,\quad d\in\N\,,\\
\label{q}
Q&=Q_0\sqcup\{b_i\colon i\to\infty\,|\, i\in I\}\,,
\\
\label{ql}
Q&=Q_0\sqcup\{b_\ell\colon\ell\to\infty\}\,, &&\ell\in I\,.
\end{align}
Quivers~\eqref{qd} and~\eqref{qld} have $md$ and $d$ added arrows, respectively; they will only be needed in Section \ref{secMultGen}. The other two quivers correspond to $d=1$. 


For $\lambda=(\lambda_0,\ldots,\lambda_{m-1})\in\C^m$ and $\alpha=(\alpha_0,\ldots,\alpha_{m-1})\in\N_0^m$ we put $\aalpha=(1,\alpha)$ and $\llambda=(-\lambda\cdot\alpha, \lambda)$ and consider representations of $\Pi^{{\llambda}}$ of dimension ${\aalpha}$.
Each representation consists of a vector space 
\begin{align}
\label{v}
 &V_\infty\oplus V_0\oplus\dots \oplus V_{m-1}\,, &&\text{where }\  V_\infty=\C\,,\quad V_i=\C^{\alpha_i}\,,\ i\in I\,,
\end{align}
together with a collection of linear maps representing the arrows of $\overline Q$, which should satisfy the relations~\eqref{relx}. It is more convenient to consider representations of $(\Pi^{\llambda})^{\mathrm{opp}}$, i.e., right $\Pi^{\llambda}$-modules (this is in agreement with the conventions in~\cite{BCE}). Therefore, analogously to~\eqref{qvar}--\eqref{gl}, we define the Calogero--Moser space associated to a framing $Q$ of $Q_0$ as
\begin{align}
\label{cms}
 &\CM_{\alpha, \lambda}(Q)=\mathtt{Rep}\big((\Pi^{\llambda})^{\mathrm{opp}}, \aalpha\big)//\mathtt{GL}(\alpha)\,, &&\alpha\in\N_0^m\,,\ \lambda\in\C^m\,. 
\end{align}
Since representations of $(\Pi^{{\llambda}})^{\mathrm{opp}}$ can be realised as dual to those of $\Pi^{\llambda}$, all the properties of $\mathtt{Rep}(\Pi^{{\llambda}}, \aalpha)$ discussed earlier remain valid in this setting. In particular, for regular $\lambda$ the varieties $\CM_{\alpha, \lambda}$ are smooth. Recall that, for a cyclic quiver $Q_0$, the root system $\Delta(Q_0)\subset\Z^{m}$ is isomorphic to the affine root system of type $\widetilde{A}_{m-1}$. Under this isomorphism, the roots in $\Delta(Q_0)$ are obtained by writing the roots of $\widetilde{A}_{m-1}$ in the basis of simple roots $\varepsilon_0,\dots, \varepsilon_{m-1}$. In particular, imaginary roots $n\delta=n(\varepsilon_0+\dots+\varepsilon_{m-1})$ correspond to $\alpha=(n,\dots, n)\in\Delta(Q_0)$. The real roots are of the form $n\delta\pm (\varepsilon_i+\dots+\varepsilon_{j-1})$ with $1\le i<j\le m-1$. Therefore, $\lambda\in\C^m$ is regular iff
\begin{align}
 &\lambda_0+\dots+\lambda_{m-1}\ne 0\,, &&n(\lambda_0+\dots+\lambda_{m-1})\ne \lambda_i+\dots+\lambda_{j-1} \label{lambdareg}
\end{align}
for all $n\in\Z$ and $1\le i<j\le m-1$. Assuming that $\lambda_i$ are the same as in~\eqref{cr1}, we conclude that the parameters $(\tau, k)\in\C^{m+1}$ are regular iff $\tau\ne 0$ and
\begin{align}
\label{kreg}
 &\frac{k_i-k_j}{\tau}\notin \frac{j-i}{m}+\Z &&\text{for}\quad 0\le i\ne j\le m-1\,. 
\end{align}
Below we will assume that the parameters are regular, unless specified otherwise.

Let us describe the varieties $\CM_{\alpha, \lambda}$ more concretely. We start by taking $Q$ as in~\eqref{q}. A right action of $\Pi^{{\llambda}}$ on a vector space~\eqref{v} consists of a collection of maps 
\begin{align*}
 &X_{a_i}\colon V_i\to V_{i+1}\,, &&X_{a_i^*}\colon V_{i+1}\to V_i\,,&&X_{b_i}\colon\C\to V_i\,, &&X_{b_i^*}\colon V_i\to \C\,.    
\end{align*}
Note that the maps go in the opposite direction to arrows because we are considering representations of the opposite algebra. To simplify the notation, put
\begin{align}
\label{maps}
 &X_i:=X_{a_i}\,, &&Y_i:=X_{a_i^*}\,, &&v_i:=X_{b_i}\,, &&w_i:=X_{b_i^*}\,.
\end{align} 
Relations~\eqref{rel} then read as follows:
\begin{align}
\label{cm1}
Y_{i}X_{i}-X_{i-1}Y_{i-1} - v_{i}w_{i}&=\lambda_i\id_{V_i} &&(i\in\mathbb Z/m\mathbb Z), \\
\label{cm2}
\sum_{i\in I}w_{i}v_{i}&=\lambda_\infty\,.
\end{align}
Note that since we are assuming that $\lambda_\infty=-\lambda\cdot\alpha$, the last relation follows from~\eqref{cm1} automatically by adding up and taking traces. Recall that the points of $\CM_{\alpha, \lambda}$ are the equivalence classes of the above data $X_i, Y_i, v_i, w_i$ under the action of $\mathtt{GL}(\alpha)$.

It will be convenient to put $V=\bigoplus_{i\in I}V_i$ and denote the projection onto the $i$-th direct summand in $V$ as $e_i$, so $e_i\big|_{V_j}=\delta_{i,j}\id_{V_i}$. Let us extend maps~\eqref{maps} up to $X_i,Y_i\colon V\to V$, $v_i\colon\C\to V$, $w_i\colon V\to\C$ in an obvious way (e.g., $Y_i\big|_{V_j}=0$ for $j\ne i+1$) and put  
\begin{align}
\label{xy}
 &X=\sum_{i\in I} X_{i}\,, &&Y=\sum_{i\in I} Y_{i}\,, &&v=\sum_{i\in I} v_i\,,&&w=\sum_{i\in I}w_i\,.
\end{align}
Then~\eqref{cm1} is equivalent to
\begin{align}
\label{mom}
YX-XY-\sum_{i\in I}v_iw_i=\sum_{i\in I}\lambda_ie_i\,,
\end{align}
and we have
\begin{align}
\label{mom1}
 &e_{i+1}X=Xe_{i}=X_i\,,\quad e_{i}Y=Ye_{i+1}=Y_i\,, &&e_iv=v_i\,,\quad we_i=w_i\,.
\end{align}
In this notation, the action of $g\in\mathtt{GL}(\alpha)\subset \End(V)$ sends $X, Y, v, w$ to
\begin{align}
\label{act}
gXg^{-1},\quad gYg^{-1},\quad gv,\quad wg^{-1}\,.
\end{align}
This action is free when $\lambda$ is regular, and its orbits are identified with the points of the Calogero--Moser space $\CM_{\alpha,\lambda}$, whose dimension of which is given by Prop.~\ref{dim}. In our case the set of vertices of $Q$ is $I\sqcup\{\infty\}$ and the dimension vector is $\aalpha=(1,\alpha)$. We have
\begin{align*}
p(\aalpha)=1+\sum_{i\in I}\alpha_i\alpha_{i+1}+\sum_{i\in I}\alpha_i-1-\sum_{i\in I}\alpha_i^2=\sum_{i\in I}\alpha_i-\frac12\sum_{i\in I}(\alpha_i-\alpha_{i+1})^2\,.
\end{align*}
Therefore, 
\begin{align}\label{dimq}
\dim \CM_{\alpha, \lambda}=2\sum_{i\in I}\alpha_i-\sum_{i\in I}(\alpha_i-\alpha_{i+1})^2\,.
\end{align}

In the case of quiver~\eqref{ql} everything looks similar except that we only have one pair $v_\ell, w_\ell$, therefore we may just put $v_i=w_i=0$ for $i\ne\ell$. This turns equations~\eqref{cm1}--\eqref{cm2} into
\begin{align}
\label{cml1}
Y_{i}X_{i}-X_{i-1}Y_{i-1} - \delta_{i,\ell}v_{\ell}w_{\ell}&=\lambda_i\id_{V_i} &&(i\in\mathbb Z/m\mathbb Z), \\
\label{cml2}
w_{\ell}v_{\ell}&=\lambda_\infty\,.
\end{align}
The relation~\eqref{mom} becomes
\begin{align}
\label{moml}
YX-XY-v_\ell w_\ell=\sum_{i\in I}\lambda_ie_i
\end{align}
and~\eqref{mom1},~\eqref{act} remain the same but with $v_i=w_i=0$ for $i\ne\ell$. Since in this case $Q$ has only one arrow going to $\infty$, we find that
\begin{align}\label{dimql}
\dim \CM_{\alpha, \lambda}=2\alpha_\ell-\sum_{i\in I}(\alpha_i-\alpha_{i+1})^2\,,
\end{align}
assuming that $(1,\alpha)$ is a positive root for $Q$ and $\lambda$ is regular.

\section{Rational solutions via Calogero--Moser spaces}
\label{secRatSol}

In this section we construct rational solutions of the hierarchy, parametrised by points of the suitable Calogero--Moser spaces. We start by setting up the notation and doing some preliminary calculations. 

\subsection{Preliminary calculations}

Let $Q$ be the quiver~\eqref{q} and $\CM_{\alpha,\lambda}$ be the space~\eqref{cms}. Take a point of $\CM_{\alpha,\lambda}$; it is represented by (isomorphism class of) data $X,Y, v, w$ satisfying~\eqref{mom}--\eqref{mom1}. Put 
\begin{align}
\label{m}
M=1-\sum_{i, j\in I}\epsilon_i w_i (X-x\id_V)^{-1}(Y-y\id_V)^{-1} v_j\epsilon_j\,.
\end{align}
Here we regard $(X-x\id_V)^{-1}$ and $(Y-y\id_V)^{-1}$ as elements of the tensor product $\mathcal P\otimes \End\, V$ by expanding $(y\id_V-Y)^{-1}$ into $\sum_{k\ge 0} Y^ky^{-k-1}$. Similarly, $e_i$, $\epsilon_i$ are identified with the elements $1\otimes e_i$ and $\epsilon_i\otimes\id_V$, so that $e_i\epsilon_i=\epsilon_ie_i=\epsilon_i\otimes e_i$, while $v_i$ and $w_i$ become maps between $\mathcal P\otimes V$ and $\mathcal P$. Therefore $M$ is an element of $\mathcal P$. It is clear that $M$ is invariant under transformations~\eqref{act}. Thus,~\eqref{m} defines a map $\CM_{\alpha, \lambda}\to 1+ \mathcal P_{-}$. 

\begin{lemma} \label{Lem-m-1}
The inverse of $M$ is given by
\begin{align}
\label{m-1}
M^{-1}=1+\sum_{i,j\in I}\epsilon_i w_i(Y-y\id_V)^{-1}(X-x\id_V)^{-1} v_j\epsilon_j\,.
\end{align}
\end{lemma} 

\noindent{\bf Proof}.
We need to check that $M^{-1}M=1$. 
Introduce the following shorthand notation: 
\begin{align*}
\wx=X-x\id_V\,,\quad \wy=Y-y\id_V\,.
\end{align*}
It easily follows from~\eqref{cr1},~\eqref{mom} that
\begin{align}
\sum_{i\in I} v_iw_i=\wy\wx- \wx\wy+\sum_{i\in I}\lambda_i(\epsilon_i-e_i)\,. \label{ijXY}
\end{align}
Next, note the following identities:
\begin{align}
\label{id1}
&e_a\epsilon_b\wx^{-1}\wy^{-1}e_c\epsilon_d=0\,,
 \qquad\quad\text{if}\  a-c\ne b-d\,,\\
\label{id2} 
&e_a\epsilon_a\wx^{-1}\wy^{-1}e_c\epsilon_c=\epsilon_ae_a\wx^{-1}\wy^{-1}e_c\,,\quad e_c\epsilon_c\wx^{-1}\wy^{-1}e_a\epsilon_a=e_c\wx^{-1}\wy^{-1}e_a\epsilon_a\,, \\
\label{id3}
&e_a\epsilon_a\wx^{-1}\wy^{-1}(\epsilon_c-e_c)=0\,,\qquad\qquad\quad (\epsilon_c-e_c)\wx^{-1}\wy^{-1}e_a\epsilon_a=0\,. 
\end{align}
Here $a,b,c,d\in I=\Z/m\Z$. To prove~\eqref{id1}, note that
\begin{multline*}
e_a\epsilon_b\wx^{-1}\wy^{-1}e_c\epsilon_d=\sum_{i,j\ge 0}e_a\epsilon_b X^iY^jx^{-1-i}y^{-1-j}e_c\epsilon_d\\= \sum_{i,j\ge 0}e_a\epsilon_b e_{c+i-j}\epsilon_{d+i-j}X^iY^jx^{-1-i}y^{-1-j}\,,
\end{multline*}
which is obviously zero if $a-c\ne b-d$ (here we used~\eqref{epsilonxy} and~\eqref{mom1}). Formula~\eqref{id2} follows from~\eqref{id1} by substituting $e_c=\sum_{i\in I} e_c\epsilon_i$ to the right hand side of~\eqref{id2}. Finally, to derive~\eqref{id3} one replaces $\epsilon_c-e_c$ by $\sum_{i\in I} (e_i\epsilon_c-e_c\epsilon_i)$ and uses~\eqref{id1}.  

Below we will repeatedly use the following corollary of~\eqref{id2}:
\begin{align}
\label{id4}
\epsilon_kw_k\wx^{-1}\wy^{-1}v_j\epsilon_j=w_k\wx^{-1}\wy^{-1}v_j\epsilon_j=\epsilon_kw_k\wx^{-1}\wy^{-1}v_j\,.
\end{align}
In other words, one can remove either $\epsilon_k$ or $\epsilon_j$ in this expression.

Now, using the formulas for $M$, $M^{-1}$ we get:
\begin{multline}
\label{cal1}
M^{-1}M=1-\sum_{i,j}\epsilon_i w_i\widehat X^{-1}\widehat Y^{-1}v_j\epsilon_j+
\sum_{i,j}\epsilon_i w_i\widehat Y^{-1}\widehat X^{-1}v_j\epsilon_j\\-\sum_{i,j,k}\epsilon_i w_i\widehat Y^{-1}\widehat X^{-1}v_k\epsilon_k w_k \widehat X^{-1}\widehat Y^{-1}v_j\epsilon_j\,.
\end{multline}
We can use~\eqref{id4} and replace the last sum by
\begin{align}
\label{cal2}
\sum_{i,j}\epsilon_i w_i\widehat Y^{-1}\widehat X^{-1}\Big(\sum_{k}v_kw_k\Big) \widehat X^{-1}\widehat Y^{-1}v_j\epsilon_j\,.
\end{align}
Replacing $\sum_k v_kw_k$ by~\eqref{ijXY} and substituting back in~\eqref{cal1}, we obtain 
\begin{align*}
M^{-1}M=1-\sum_{i,j,k\in I}\lambda_k \epsilon_i w_i\widehat Y^{-1}\widehat X^{-1}(\epsilon_k-e_k)\widehat X^{-1}\widehat Y^{-1}v_j\epsilon_j\,.
\end{align*}
By~\eqref{id3},
\[
(\epsilon_k-e_k)\widehat X^{-1}\widehat Y^{-1}v_j\epsilon_j=(\epsilon_k-e_k)\widehat X^{-1}\widehat Y^{-1}e_j\epsilon_jv_j=0\,.
\] 
We conclude that $M^{-1}M=1$, as needed. \qed

\begin{lemma}\label{L-} Let $L=MyM^{-1}$. For $k\in\N$, we have
\begin{multline}
 (L^k)_-=\sum_{i,j}\epsilon_{i}w_{i-k}Y^k\wy^{-1}\wx^{-1}v_j\epsilon_j
-\sum_{i,j}\epsilon_iw_i\wx^{-1}\wy^{-1}Y^kv_{j+k}\epsilon_{j} \\
-\sum_{i,j,\ell}\epsilon_iw_i\wx^{-1}\wy^{-1}Y^k v_\ell w_{\ell-k}\wy^{-1}\wx^{-1}v_j\epsilon_j\\+\sum_{i,j,\ell}\sum_{a=0}^{k-1}\epsilon_iw_i\wx^{-1}Y^a v_\ell w_{\ell-k}Y^{k-1-a}\wy^{-1}\wx^{-1}v_j\epsilon_j\,. \label{LkminusNE}
\end{multline}
\end{lemma}

\noindent{\bf Proof.} Write $M, M^{-1}$ as
\begin{align*}
M=1-\sum_{i,j}\mu_{ij}\,,\quad M^{-1}=1+\sum_{i,j}\tilde\mu_{ij}\,,
\end{align*}
where
\begin{align*}
\mu_{ij}=\epsilon_iw_i\wx^{-1}\wy^{-1}v_j\epsilon_j\,,\quad \tilde\mu_{ij}=\epsilon_iw_i\wy^{-1}\wx^{-1}v_j\epsilon_j\,.
\end{align*}
Therefore,
\begin{align} \label{Lk}
 L^k=My^kM^{-1}=y^k+\sum_{i, j}y^k\tilde\mu_{ij}-\sum_{i, j}\mu_{ij}y^k-\sum_{iji'j'}\mu_{ij}y^k\tilde\mu_{i'j'}\,.
\end{align}
Taking into account that $\epsilon_jy^k=y^k\epsilon_{j-k}$, one gets 
\begin{align}
\label{mu1}
 \sum_{i, j}y^k\tilde\mu_{ij}&=\sum_{i,j}\epsilon_{i}w_{i-k}y^k\wy^{-1}\wx^{-1}v_j\epsilon_j\,, \\
\label{mu2}
-\sum_{i,j}\mu_{ij}y^k&=-\sum_{i,j}\epsilon_iw_i\wx^{-1}\wy^{-1}y^kv_{j+k}\epsilon_{j}\,,\\
\label{mu3}
-\sum_{i, j, i', j'}\mu_{ij}y^k\tilde\mu_{i'j'}&=-\sum_{i, j, \ell}\epsilon_iw_i\wx^{-1}\wy^{-1}\epsilon_\ell v_\ell y^kw_{\ell-k}\wy^{-1}\wx^{-1}v_j\epsilon_j\,.
\end{align}
By~\eqref{id4}, we can remove $\epsilon_\ell$ in the last formula to get
\begin{align*}
-\sum_{i, j, i', j'}\mu_{ij}y^k\tilde\mu_{i'j'}=-\sum_{i, j, \ell}\epsilon_iw_i\wx^{-1}\wy^{-1}y^kv_\ell w_{\ell-k}\wy^{-1}\wx^{-1}v_j\epsilon_j\,.
\end{align*}
From an elementary formula $Y^k-y^k\id_V=(Y-y\id_V)\sum_{a=0}^{k-1}y^{k-1-a}Y^a$
we obtain:
\begin{equation}\label{4.13a}
y^k\wy^{-1}=\wy^{-1}y^k=\wy^{-1}Y^k-\sum_{a=0}^{k-1}Y^ay^{k-1-a}\,,\qquad (y^k\wy^{-1})_-=Y^k\wy^{-1}\,.
\end{equation}
Similarly,
\begin{align*}
\wy^{-1}y^kv_\ell w_{\ell-k}\wy^{-1}=\wy^{-1}Y^kv_\ell w_{\ell-k}\wy^{-1}- \sum_{a=0}^{k-1}Y^av_\ell w_{\ell-k}y^{k-1-a}\wy^{-1}\,,
\end{align*}
and by applying the formula \eqref{4.13a} we obtain
\begin{align*}
(\wy^{-1}y^kv_\ell w_{\ell-k}\wy^{-1})_-=\wy^{-1}Y^kv_\ell w_{\ell-k}\wy^{-1}- \sum_{a=0}^{k-1}Y^av_\ell w_{\ell-k}Y^{k-1-a}\wy^{-1}\,.
\end{align*}
Using these facts, we compute the negative part of~\eqref{mu1}--\eqref{mu3} which results in~\eqref{LkminusNE}. \qed

\subsection{Solutions for the hierarchy}

Let $Q$ be the quiver~\eqref{q} and $\aalpha=(1, \alpha)$ with $\alpha\in \N_0^I$. Consider $\Rep(\C\overline Q^{\mathrm{opp}},\aalpha)$; a point of $\Rep(\C\overline Q^{\mathrm{opp}},\aalpha)$ is a collection of maps 
\begin{align*}
 &X_i\colon V_i\to V_{i+1}\,, &&Y_i\colon V_{i+1}\to V_i\,, &&v_i\colon \C\to V_i\,, &&w_i\colon V_i\to \C\,.  
\end{align*}
Thus, this is an affine space of dimension $2\sum_{i\in I}(\alpha_i\alpha_{i+1}+\alpha_i)$. It has a natural symplectic structure~\eqref{symp}, 
\begin{align}
\label{symp1}
\omega=\sum_{i\in I}\left(\tr(\mathrm{d}Y_i \wedge \mathrm{d}X_{i})+\mathrm{d}w_i \wedge \mathrm{d}v_i\right)\,.
\end{align}
The corresponding Poisson bracket is given by
\begin{align}
\label{poisson}
 &\big\{(X_i)_{ab}, (Y_j)_{cd}\big\}=\delta_{ij}\delta_{ad}\delta_{bc}\,, &&\big\{(v_i)_{a}, (w_j)_{b}\big\}=\delta_{ij}\delta_{ab}\,.
\end{align}
(All remaining brackets are zero.)

As before, we will combine $X_i$, $Y_j$, $v_i$, $w_i$ into $X, Y\in \End V$ and $v\colon \C\to V$, $w\colon \C\to V$, see~\eqref{xy}. 
For $k\in \N$, define a flow on $\Rep(\C\overline Q^{\mathrm{opp}},\aalpha)$ by the formulas
\begin{align}
  \dot{X} &=-\sum_{\ell\in I}\sum_{a=0}^{k-1} Y^av_\ell w_{\ell-k}Y^{k-a-1}\,, 
\label{sFlowXY} 
 \qquad\dot{Y}=0\,,\\
 \dot{v}&=-Y^kv\,,\qquad\dot{w}=wY^k\,, \label{sFlowVW}
\end{align}
where the dot denotes the derivative with respect to the time variable $t_k$. Equivalently, 
\begin{align}
 \dot{X_i}&=-\sum_{a=0}^{k-1} Y_{i+1}Y_{i+2}\cdots Y_{i+a} v_{i+a+1} w_{i+a+1-k}Y_{i+a+1-k}\cdots Y_{i-1}\,, \label{sFlowXi} \\
 \dot{Y_i}&=0\,, \label{sFlowYi} \\
 \dot{v}_{i}&=-Y^kv_{i+k}=-Y_{i}Y_{i+1}\cdots Y_{i+k-1} v_{i+k}\,,\label{sFlowVi} \\
 \dot{w}_{i}&=w_{i-k}Y^k=w_{i-k}Y_{i-k}Y_{i-k+1}\cdots Y_{i-1}\,. \label{sFlowWi}
\end{align}
It is easy to check that this is a Hamiltonian flow, with the Hamiltonian
\begin{align}
\label{ham}
H_k=-wY^kv=-\sum_{i\in I} w_iY_i\dots Y_{i+k-1}v_{i+k}\,.
\end{align}
The Hamiltonian $H_k$ is invariant under the action~\eqref{act} of $G(\aalpha)$, therefore, the moment map is a first integral, thus the flow descends to a Hamiltonian flow on a symplectic quotient $\CM_{\alpha, \lambda}$. (This can also be checked directly by differentiating the relations~\eqref{cm1}--\eqref{cm2}.)

\begin{prop} \label{prop4.3}
Let $L=MyM^{-1}$ where $M$, $M^{-1}$ are as in~\eqref{m} and~\eqref{m-1}. If $X, Y, v, w$ satisfy the equations~\eqref{sFlowXY}--\eqref{sFlowVW}, then $L$ satisfies $\dot{L}=[(L^k)_+, L]$.
\end{prop}

\noindent{\bf Proof.} According to Prop.~\ref{wave}, it suffices to check that
\begin{align}
\label{meq}
 \dot{M}=-(L^k)_- M.
\end{align}

Let us start by calculating the right-hand side. Note that the second and third sum in~\eqref{LkminusNE} are nothing but 
\begin{align*}
\Big(-\sum_{i,j}\epsilon_iw_i\wx^{-1}\wy^{-1}Y^kv_{j+k}\epsilon_{j} \Big)M^{-1}\,.
\end{align*}
Therefore, substituting them in the right-hand side of~\eqref{meq} results in
\begin{align}
\label{terms1}
\sum_{i,j}\epsilon_iw_i\wx^{-1}\wy^{-1}Y^kv_{j+k}\epsilon_{j}\,.
\end{align}
Multiplying the remaining terms in~\eqref{LkminusNE} by $M$ gives:
\begin{multline}
\label{terms2}
\Big(\sum_{i,j}\epsilon_{i}w_{i-k}Y^k\wy^{-1}\wx^{-1}v_j\epsilon_j
+\sum_{i,j,\ell}\sum_{a=0}^{k-1}\epsilon_iw_i\wx^{-1}Y^a v_\ell w_{\ell-k}Y^{k-1-a}\wy^{-1}\wx^{-1}v_j\epsilon_j\Big)M\\
=\sum_{i,j}\epsilon_{i}w_{i-k}Y^k\wy^{-1}\wx^{-1}v_j\epsilon_j
+\sum_{i,j,\ell}\sum_{a=0}^{k-1}\epsilon_iw_i\wx^{-1}Y^a v_\ell w_{\ell-k}Y^{k-1-a}\wy^{-1}\wx^{-1}v_j\epsilon_j\\
-\sum_{i,j, r}\epsilon_{i}w_{i-k}Y^k\wy^{-1}\wx^{-1}v_r w_r\wx^{-1}\wy^{-1}v_j\epsilon_j\\
-\sum_{i,j,\ell, r}\sum_{a=0}^{k-1}\epsilon_iw_i\wx^{-1}Y^a v_\ell w_{\ell-k}Y^{k-1-a}\wy^{-1}\wx^{-1}v_rw_r\wx^{-1}\wy^{-1}v_j\epsilon_j\,.
\end{multline}
Replacing $\sum_r v_r w_r$ in the last two terms by~\eqref{ijXY} and applying~\eqref{id3} and~\eqref{id4}, we are left after cancellations with the following terms:
\begin{align}
\label{terms3}
\sum_{i,j}\epsilon_{i}w_{i-k}Y^k\wx^{-1}\wy^{-1}v_j\epsilon_j
+\sum_{i,j,\ell}\sum_{a=0}^{k-1}\epsilon_iw_i\wx^{-1}Y^a v_\ell w_{\ell-k}Y^{k-1-a}\wx^{-1}\wy^{-1}v_j\epsilon_j\,.
\end{align}
Now, the right-hand side of~\eqref{meq} is obtained by subtracting~\eqref{terms3} from~\eqref{terms1}, so
\begin{align}
-(L^k)_-M=&\sum_{i,j}\epsilon_iw_i\wx^{-1}\wy^{-1}Y^kv_{j+k}\epsilon_{j}-\sum_{i,j}\epsilon_{i}w_{i-k}Y^k\wx^{-1}\wy^{-1}v_j\epsilon_j\nonumber\\\label{rhs}
&-\sum_{i,j,\ell}\sum_{a=0}^{k-1}\epsilon_iw_i\wx^{-1}Y^a v_\ell w_{\ell-k}Y^{k-1-a}\wx^{-1}\wy^{-1}v_j\epsilon_j\,.
\end{align}
On the other hand, from the formula for $M$ we get that
\begin{align*}
 \dot{M}=-\sum_{i, j}\epsilon_i{w_i}\wx^{-1}\wy^{-1}\dot{v}_{j}\epsilon_j-\sum_{i, j}\epsilon_i\dot{w}_{i}\wx^{-1}\wy^{-1}v_j\epsilon_j+\sum_{i,j} \epsilon_i{w_i}\wx^{-1}\dot{X}\wx^{-1}\wy^{-1}{v_j}\epsilon_j\,.
\end{align*}
Now a quick look at~\eqref{rhs} and~\eqref{sFlowXY},~\eqref{sFlowVi},~\eqref{sFlowWi} confirms that $\dot{M}=-(L^k)_-M$, as needed.
\qed

With the flow~\eqref{sFlowXY}--\eqref{sFlowVW} on $\CM_{\alpha, \lambda}$ corresponding to the $k$-th flow of the hierarchy, we obtain a commuting family of Hamiltonian flows on $\CM_{\alpha, \lambda}$ labelled by $k$. In fact, one can confirm by a direct calculation that the Hamiltonians~\eqref{ham} commute, that is, 
\begin{align}
\label{comH}  
\{H_j, H_k\}=0\qquad\text{for all}\ \,j, k\,.
\end{align}
Thus, we arrive at the following result.

\begin{theorem}\label{sol}
Under the map $\CM_{\alpha, \lambda}\to\mathcal P$ which sends $(X,Y,v,w)$ to $L=MyM^{-1}$, the commuting flows with the Hamiltonians $H_k$, $k=1,2,\dots$, get identified with the flows of the hierarchy~\eqref{kp}.  
\end{theorem}
\qed

Let us integrate the flows~\eqref{sFlowXY}--\eqref{sFlowVW} simultaneously for all $k$. First we integrate~\eqref{sFlowVW} to find that 
\begin{align}
v=e^{-\sum_{k}t_kY^k}v(0)\,,\quad w=w(0)e^{\sum_{k}t_kY^k}\,,\quad vw=e^{-\sum_{k} t_k\ad_{Y^k}}\big(v(0)w(0)\big)\,, \label{vwsol}
\end{align}
where $Y=Y(0)$ is constant and we use the standard notation for $\ad_A B=AB-BA$. Next, note that~\eqref{sFlowXY} can be written as
\begin{align*}
\frac{\partial X}{\partial t_k}=-\sum_{i\in I} e_{i+1}\Big(\sum_{a=0}^{k-1}Y^avwY^{k-1-a}\Big)e_{i}= -\sum_{i\in I}e_{i+1}(\ad_Y)^{-1}\ad_{Y^k}(vw)e_i\,.
\end{align*}
Substituting the expression for $vw$ and integrating, we find that
\begin{align}
 X=X(0)+\sum_{i\in I} e_{i+1}\left\{(\ad_Y)^{-1}\big(e^{-\sum_{k} t_k\ad_{Y^k}}-1\big)\big(v(0)w(0)\big)\right\}e_i\,. \label{Xsol}
\end{align}
Note that $(\ad_Y)^{-1}(e^{-\sum_{k} t_k\ad_{Y^k}}-1)$ is well defined: after expanding the exponential into a series each term contains a well-defined factor $(\ad_Y)^{-1}\ad_{Y^k}$.  

Let us remark on the flows with $k=mj$, $j\in \N$. In this case 
\begin{align*}
 H_{mj}=-wY^{mj}v=-\sum_{i\in I} w_iY^{mj}v_i=-\tr\big(Y^{mj}\sum_i v_iw_i\big)\,.
\end{align*} 
On $\CM_{\alpha, \lambda}$~\eqref{mom} holds, so $H_{mj}$ can be replaced by
\begin{equation}\label{hk}
\widetilde H_{mj}=\tr\bigg(Y^{mj}\Big(XY-YX+\sum_{i\in I} \lambda_ie_i\Big)\bigg)= \sum_{i\in I}\lambda_i\tr(e_iY^{mj}e_i)=\frac{|\lambda|}{m}\, \tr Y^{mj}\,,
\end{equation}
where $|\lambda|:=\sum_{i\in I}\lambda_i$. According to~\eqref{cr1}, $|\lambda|=-m\tau\ne0$. Therefore, the flow~\eqref{sFlowXY}--\eqref{sFlowVW} can be replaced by
 \begin{align}
  \dot{X} &=j|\lambda|Y^{mj-1}\,, \label{FlowXY} 
 & \dot{Y}&=0\,,\\
 \dot{v}&=0\,,
 & \dot{w}&=0\,. \label{FlowVW}
\end{align}
As a result, the $mj$-th flow is linear: $X(t_{mj})=X(0)+j|\lambda|t_{mj}Y^{mj-1}$.

\subsection{Solutions for the spherical sub-hierarchy}   
\label{secSolSubh}

Given $\ell\in I$, consider the quiver~\eqref{ql} and let $\CM_{\alpha,\lambda}$ be a Calogero--Moser space~\eqref{cms}. Recall that the description of $\CM_{\alpha,\lambda}$ is almost the same as for the quiver~\eqref{q}, see~\eqref{cml1}--\eqref{moml}, the only difference being that $v_i=w_i=0$ for $i\ne \ell$. Therefore, in this case $v=v_\ell$, $w=w_\ell$ and we simply need to make appropriate changes in the previous formulas. 

Hence, given a point $(X,Y,v,w)$ of $\CM_{\alpha, \lambda}$, we put
\begin{align}
\label{ml}
M=1-\epsilon_\ell w_\ell(X-x\id_V)^{-1}(Y-y\id_V)^{-1}v_\ell \epsilon_\ell\,.
\end{align}
This defines a map $\CM_{\alpha, \lambda}\to 1+ \mathcal P_{-}$. By Lemma~\ref{Lem-m-1}, the inverse of $M$ is given by
\begin{align}
\label{ml-1}
M^{-1}=1+\epsilon_\ell w_\ell(Y-y\id_V)^{-1}(X-x\id_V)^{-1}v_\ell \epsilon_\ell\,.
\end{align}
It is obvious from these formulas that $\epsilon_i M=M\epsilon_i$ and $\epsilon_i M^{-1}=M^{-1}\epsilon_i$ for $i\in I$. Hence, $L=MyM^{-1}$ and $\widetilde L=L^m=My^mM^{-1}$ will satisfy
\begin{align*}
 &\epsilon_i L = L \epsilon_{i-1}\,, &&\epsilon_i\widetilde L=\widetilde L\epsilon_i\,.
\end{align*}
Now, using $(X,Y,v,w)$ as the initial data, we can construct a solution to the spherical sub-hierarchy. As explained above, when restricting to the flows with $k=mj$, the dependence of $(X,Y,v,w)$ on $t=(t_m, t_{2m}, \dots)$ reduces to
\begin{align}\label{hd}
 &X(t)=X(0)+|\lambda| \sum_{j\ge 1} j t_{mj}Y^{mj-1}\,, &&Y=Y(0)\,, &&v=v(0)\,, && w=w(0)\,.
\end{align}
Clearly, the condition $v_i=w_i=0$ for $i\ne\ell$ remains valid for all times. As a result, the corresponding $L(t)$, $\widetilde L(t)$ and $\widetilde L_\ell(t)=\epsilon_\ell\widetilde L(t)$ will satisfy the equations~\eqref{kpm},~\eqref{kpmm} and~\eqref{kpmmm}. 

\begin{remark} 
Up to a change of notation, the element $M$ defined in \eqref{ml} is the same as $\kappa$ in \cite[Theorem 2.6]{Esh};  interpreting it as a ``pseudo-differential operator" in a suitable sense was the starting point in our work.
\end{remark}

\begin{remark} The Calogero--Moser spaces \eqref{cms} associated with the quivers \eqref{q}--\eqref{ql} are particular cases of the quiver varieties that appear in \cite{BGK1, BGK2}. Namely, they correspond to setting the $\Gamma$-module $W$ in \cite[(2.16)]{BGK2} as $W=\C\Gamma$ and $W=\C\epsilon_\ell$, respectively, where $\C\epsilon_\ell$ is viewed as a one-dimensional $\Gamma$-module. Similarly, the Calogero--Moser spaces for the quivers \eqref{qd}--\eqref{qld} correspond to $W=\C\Gamma^{\oplus d}$ and $W=(\C\epsilon_\ell)^{\oplus d}$. 
\end{remark}

\begin{remark}
Above we assumed that $\lambda$ is regular, which puts restrictions~\eqref{kreg} on parameters $k_0, \dots, k_{m-1}$ for the Cherednik algebra $H_{\tau, k}$. Let us remark on the situation in general, for arbitrary $k_i$ and $\tau\ne 0$. In this case one may have reducible representations in $\mathtt{Rep}(\Pi^{\llambda}, \aalpha)$. Non semi-simple representations are irrelevant because they are indistinguishable from semi-simple ones under the map $(X,Y,v,w)\mapsto M$. On the other hand, if $(X,Y,v,w)$ describes a semi-simple reducible $\Pi^{\llambda}$-module, then it is a direct sum of simples, exactly one of which must contain $V_\infty=\C$. The formula \eqref{m} makes it clear that $M$ is fully determined by this simple summand. Let $\aalpha'=(1, \alpha')$ be the dimension vector of this summand, so it corresponds to some point $(X',Y',v',w')\in\mathtt{Rep}(\Pi^{\llambda}, \aalpha')$. Then the two solutions of the hierarchy built from $(X,Y,v,w)$ and $(X', Y', v', w')$ coincide. Therefore, the solutions of the hierarchy are parametrised by the locus of simple representations: 
$$\bigsqcup_{\aalpha=(1, \alpha)}\Irr(\Pi^{\llambda}, \aalpha)// \mathtt{GL}(\alpha)\,.$$
For a description of the set ${\Sigma}_\lambda\subset \Delta_+(Q)$ of all $\aalpha$ such that $\Irr(\Pi^{\llambda}, \aalpha)\ne \emptyset$,  see Theorem 1.2 of~\cite{CB2}.   

\end{remark}

\begin{remark} \label{remGenQuiver}
The family of commuting Hamiltonians~\eqref{ham} can be generalised to any quiver $Q_0$ with a vertex set $I$. Namely, suppose $Q$ is obtained from $Q_0$ by adding a vertex $\infty$ and one arrow $v_i\colon \infty \to i$ for every $i\in I$. For $k\in \N$, let $h_k=\sum_{i,j\in I} v^*_ja_k \dots a_1v_i\in\C[\overline Q]$, where the sum is taken over all paths $a_k\dots a_1$ length $k$ in $Q_0$. It is a pleasant exercise to check that $[h_k, h_l]_{NL}=0$, where $[\cdot , \cdot]_{NL}$ denotes the necklace Lie bracket~\cite{BLB}. According to the general principle~\cite{VdB, CB4}, for any representation $\rho\in\mathtt{Rep}(\C\overline Q, \alpha)$ the functions $H_k=\tr \rho(h_k)$ Poisson-commute both on $\mathtt{Rep}(\C\overline Q, \alpha)$ and on the symplectic quotients $\mathtt{Rep}(\Pi^{\lambda}, \alpha)//G(\aalpha)$.  
\end{remark}

\section{Link to the Calogero--Moser system} \label{sec5}

In this section we explain how the dynamics~\eqref{sFlowXY}--\eqref{sFlowVW} is related to the classical Calogero--Moser system for the wreath-product $\Z_m\wr S_n$. For this purpose, we will make the special choice $\aalpha=(1, n\delta)$ where $\delta=(1, \dots , 1)$. Consider the varieties $\CM_{n\delta, \lambda}$ which will be abbreviated to $\CM_{n,\lambda}$ for simplicity. Everything still depends on the choice of framing, so we have two versions of $\CM_{n,\lambda}$, corresponding to quivers~\eqref{q},~\eqref{ql}. We will deal first with the case~\eqref{ql}.

Let $\Gamma$ denote, as previously, the cyclic group of order $m$ with a generator ${s}$ acting on $\C$ by a fixed primitive $m$th root of unity, $\mu$. The group $\Gamma^n$ acts naturally on $\mathfrak h:=\C^n$. The wreath product $G=\Gamma\wr S_n$ is, by definition, the semidirect product $\, S_{n} \ltimes \Gamma^n \, $, where $S_n$ acts in the obvious way on $\mathfrak h$ by permuting the coordinates. Write ${s}_{ij}$ for the elementary transpositions in $S_n$, and ${s}_i$ for a copy of ${s}\in\Gamma$ placed in the $i$th factor in $\Gamma^n$. The group $G$ is generated by all ${s}_{ij}$ and ${s}_i$. This is the group $G(m, 1, n)$ in the standard classification of complex reflection groups~\cite{ST}; in the case $m=2$ it is the Coxeter group of type $B_n$.

\subsection{Darboux coordinates on $\CM_{n,\lambda}$}
Consider the quiver~\eqref{ql}, where we set $\ell=0$ for simplicity. From the formula~\eqref{dimql}, we find that $\dim \CM_{n,\lambda}=2n$. Points of $\CM_{n,\lambda}$ are represented by linear maps
\begin{align*}
 &X_{i}\colon V_i\to V_{i+1}\,, &&Y_{i}\colon V_{i+1}\to V_i\,, &&v_{0}\colon \C\to V_0\,, &&w_{0}\colon V_0\to \C\,,    
\end{align*}
where $\dim V_i=n$. These maps should satisfy the relations~\eqref{cml1}:
\begin{align}
\label{cm0}
 &Y_{i}X_{i}-X_{i-1}Y_{i-1} - \delta_{i0}v_{0}w_{0}=\lambda_i\id_{V_i}  &&(i\in\mathbb Z/m\mathbb Z)\,.
\end{align}
To construct local coordinates on $\CM_{n,\lambda}$, we look at the subset $\CM'_{n,\lambda}\subset \CM_{n,\lambda}$ where the transformation $X_{m-1}\dots X_1X_0\in \End(V_0)$ is diagonalisable and invertible. Suppose that this transformation has eigenvalues $\nu_1^m,\ldots,\nu_n^m$. Then one can choose bases in each of $V_i$ so that all $X_i$ look identical:
\begin{equation}
\label{slice0}
X_i=\mathrm{diag}({\nu}_1, \dots, {\nu}_n)\,.
\end{equation}
Then from~\eqref{cm0} we obtain that
\begin{align}
\label{mr1}
(Y_i)_{jj}\nu_j-\nu_j(Y_{i-1})_{jj}-\delta_{i0}(v_0w_0)_{jj}&=\lambda_i\,,\\
\label{mr2}
(Y_i)_{jk}\nu_k-\nu_j(Y_{i-1})_{jk}-\delta_{i0}(v_0w_0)_{jk}&=0 &&(j\ne k)\,.
\end{align}
Summing~\eqref{mr1} over $i$ we obtain that for any $j$, 
\begin{align*}
 &(v_0w_0)_{jj}=-|\lambda|\,, &&\text{where}\quad |\lambda|=\sum_{i\in I}\lambda_i\ne 0\,.
\end{align*} 
Therefore, by rescaling the bases in each of $V_i$ we may achieve that
\begin{align}
\label{slice1}
(v_0)_j=1\,, \quad (w_0)_j=-|\lambda|\,,\quad\text{and therefore $(v_0w_0)_{jk}=-|\lambda|$ for all $j,k$}\,.
\end{align}
The matrices $Y_i$ are then easily determined from relations~\eqref{mr1}--\eqref{mr2}. Namely, we must have $\nu_j^m\ne \nu_k^m$ for $j\ne k$, and the matrices $Y_i$ can be chosen in the following form:
\begin{align}
\label{slice2}
&(Y_i)_{jj}=\mu_j+\frac{1}{\nu_j}(a+\sum_{r=0}^{i}\lambda_r)\,,\qquad a=-\sum_{r=0}^{m-1}\frac{m-r}{m}\lambda_r\,,\\
 \label{slice3}
&(Y_i)_{jk}=|\lambda|\,\frac{\nu_j^{i}\nu_k^{m-i-1}}{\nu_j^m-\nu_k^m}\qquad\qquad (j\ne k)\,,
\end{align}
where $i=0, 1, \dots, m-1$. Here $\mu_1, \dots, \mu_n$ are arbitrary, so $\{\mu_j, \nu_j\}$ are local coordinates on $\CM'_{n,\lambda}$.

Note that the set of eigenvalues $\{\nu_1^m,\ldots,\nu_n^m\}$ determines $\nu_1,\ldots,\nu_n$ only up to permutations and multiplication by $m$th roots of unity. As a result, the coordinates $\nu_j$, $\mu_j$ on $\CM'_{n,\lambda}$ are defined up to the action of $G=S_n\ltimes\Gamma^n\subset\GL(n,\C)$ that simultaneously permutes $\{\nu_j\}$, $\{\mu_j\}$ and replaces $\nu_j$, $\mu_j$ by $\mu^{r_j}\nu_j$, $\mu^{-r_j}\mu_j$. This action comes from the $G$-action on $\Rep\big((\Pi^\llambda)^{\mathrm{opp}},\aalpha\big)$ via the embedding $G\hookrightarrow\GL(\alpha)=\GL(n,\C)^m$ given by
\begin{align} \label{embed}
 &\sigma\mapsto(\sigma,\ldots,\sigma)\,, &&s_1^{r_1}\cdots s_n^{r_n}\mapsto(1,F,\ldots,F^{m-1}),
\end{align}
where $\sigma\in S_n$ and $F=\diag(\mu^{r_1},\ldots,\mu^{r_n})$.

It follows that $\CM'_{n,\lambda}$ is isomorphic to $\hreg \times \mathfrak h^*/G$, where
\begin{equation} \label{hreg}
 \hreg:=\big\{x\in\mathfrak h\,\mid\, \prod_{i=1}^nx_i\prod_{i<j}(x_i^m-x_j^m)\ne 0\big\}\,.
\end{equation} 
 
Recall that $\CM_{n,\lambda}$ inherits a symplectic structure from the standard symplectic structure on the space $\mathtt{Rep}\big(\C\overline{Q}^{\mathrm{opp}}, \aalpha\big)$,
\begin{equation}
\label{symp2}
\omega=\sum_{i\in I}\tr(\mathrm{d}Y_i \wedge \mathrm{d}X_{i})+\mathrm{d}w_0 \wedge \mathrm{d}v_0\,.
\end{equation}
Denote this symplectic structure as $\omega_{\CM_{n,\lambda}}$. 
\begin{prop}\label{sf}
The above coordinates $\{\mu_i, \nu_i\}$ are Darboux coordinates in the sense that 
\begin{align*}
 \omega_{\CM_{n,\lambda}}=m\sum_{i=1}^n \mathrm{d}\mu_i \wedge \mathrm{d}\nu_i\,. 
\end{align*} 
\end{prop}

\noindent{\bf Proof.} Recall that $\CM_{n,\lambda}$ is obtained from $\Rep(\C Q^{\mathrm{opp}}, \aalpha)$ by Hamiltonian reduction, see \eqref{sympq}. Using the coordinates $\{\mu_i, \nu_i\}$ we identify $\CM'_{n,\lambda}$ locally with an open subset in $\hreg\times\mathfrak h^*$. 
Consider the mapping $\phi\colon\hreg\times\mathfrak h^*\to\mathtt{Rep}(\C\overline{Q}^{\mathrm{opp}}, \aalpha)$, defined by the formulas~\eqref{slice0} and~\eqref{slice1}--\eqref{slice3}. Locally, this represents a transversal slice for the action of $\GL(n\delta)=\GL(n, \C)^m$ on $\Rep\big((\Pi^{\llambda})^{\mathrm{opp}}, \aalpha\big)\subset\mathtt{Rep}(\C\overline{Q}^{\mathrm{opp}}, \aalpha)$. Therefore, the symplectic structure on $\CM'_{n,\lambda}$ is the pullback $\phi^*\omega$ of the form~\eqref{symp2}. To calculate it, we substitute the expressions for $X_i, Y_i, v_0, w_0$ into~\eqref{symp1}. With $v_0$ and $w_0$ being constant and $\partial X_i/\partial{\mu}_i=0$, 
we get
\begin{equation*}
\phi^*\omega=\sum_{i\in I}\sum_{j<k}\tr\left(\frac{\partial Y_i}{\partial\nu_j}\frac{\partial X_i}{\partial\nu_k}-\frac{\partial Y_i}{\partial\nu_k}\frac{\partial X_i}{\partial\nu_j}\right)\mathrm{d}\nu_j \wedge \mathrm{d}\nu_{k}+
\sum_{i\in I}\sum_{j,k}\tr\left(\frac{\partial Y_i}{\partial\mu_j}\frac{\partial X_i}{\partial\nu_k}\right)\mathrm{d}\mu_j \wedge \mathrm{d}\nu_k\,.
\end{equation*}
Since each $X_i$ is a diagonal matrix, the off-diagonal part of $Y_i$ does not play any role and can be ignored. Then it is easy to see that the only non-zero terms are those with $j=k$, where we have $\tr\left(\frac{\partial Y_i}{\partial\mu_j}\frac{\partial X_i}{\partial\nu_j}\right)=1$. This leads to the result of the proposition. \qed

\begin{remark}\label{shi}
 The symplectic form $\omega_{\CM_{n,\lambda}}$ does not change under the transformation $\mu_j\to\mu_j-a\nu_j^{-1}$, therefore the parameter $a$ in the formula for $(Y_i)_{jj}$ can be eliminated. However, the choice of $a$ in~\eqref{slice2} ensures that the Hamiltonians $\frac1m\tr Y^{mp}$ have the form $\sum_{j=1}^n\mu_j^{pm}$ plus the terms of degree $\le mp-2$ in momenta~$\mu_j$.
\end{remark}

\subsection{Cherednik algebra and Calogero--Moser spaces}

In this section we recall the link between the Calogero--Moser spaces and Cherednik algebra, following \cite{EG, EGGO}, see also \cite{BCE}. Our notation will mostly follow that of \cite{BCE}. 

Consider the group $G=S_n\ltimes\Gamma^n$ acting on $\mathfrak h=\C^n$ and, a by natural extension, on the tensor algebra $T(\mathfrak h\oplus\mathfrak h^*)$. 
Fix parameters $t, c_{00}, c_1, \dots, c_{m-1}\in\C$, and let $c:=(c_{00},c_1, \dots, c_{m-1})$. The Cherednik algebra $H_{t,c}(G)$ is then defined as the quotient of $T(\mathfrak h\oplus\mathfrak h^*)* G$ by the following relations:
\begin{align}
 &[x_i,x_j]=0, \qquad [y_i,y_j]=0,\qquad 1\le i, j\le n\,,\nonumber \\
 &[y_i,x_i]=t+c_{00} \sum_{j\ne i}\sum_{l=0}^{m-1}{s}_{ij}{s}_i^{l}{s}_j^{-l}
+\sum_{l=1}^{m-1}c_l{s}_i^l, \qquad 1\le i\le n\,,\label{xiyi}\\
 &[y_i,x_j]=-c_{00} \sum_{l=0}^{m-1}{s}_{ij}\mu^{l}{s}_i^{l}{s}_j^{-l} \,,\qquad 1\le i\ne j\le n\,. \nonumber
\end{align}
Here $\{x_i\}$ and $\{y_i\}$ are the standard dual bases of $\mathfrak h^*$ and $\mathfrak h$, respectively. Compared to the notation of \cite{BCE}, our $\mu$, $c_{00}$ and $s_i$ correspond to $\epsilon^{-1}$, $k$ and  $\alpha_i$ in {\it loc. cit.} Note that there is a misprint in the formulas (34) in \cite{BCE}: the correct expression for $[y_i, x_j]$ should have $\epsilon^{-l}$ in place of $\epsilon^l$. Note also that $c_{00}$ in \eqref{xiyi} is denoted by $k/2$ in \cite{EGGO}.

Next, let us describe following \cite{EGGO} a functorial link between representations of $H_{t,c}$ and a suitable deformed preprojective algebra $\Pi^{\llambda}$. To this end, let $G'=S_{n-1}\ltimes\Gamma^{n-1}$ be the subgroup of $G$ which stabilises the first basis vector in $\mathfrak h=\C^n$. Write $ \bold{e} := \frac{1}{ | G | } \sum_{g \in G} g $ and $ \bold{e'} := \frac{1}{ | G' | } \sum_{g \in G'} g $ for the symmetrising idempotents in $ \mathbb C G$ and $\C G'$, respectively. Now, recall the idempotents \eqref{id} and, 
for each $i=0,1,\ldots, m-1$, define $\ee_{i}=(\epsilon_i)_1\,\bold{e'}$, where $(\epsilon_i)_1$ is a copy of the idempotent $\epsilon_i\in\C\Gamma$ placed in the first component of $\C\Gamma^n$.  
One can check easily that $ \ee_i $ are idempotents in $ \C G\subset H_{t,c}$
satisfying the relations
\begin{equation*}
 x_{1} \cdot \ee_{i+1}=
\ee_{i} \cdot x_{1}\, , \quad  y_{1} \cdot \ee_{i} = \ee_{i+1}\cdot
y_{1}\ .
\end{equation*}
Now, take an arbitrary (left) $H_{t,c}$-module $E$ and put $V_\infty=\bold{e} E$ and $V_i=\ee_{-i} E$ for $i=0, \dots, m-1$. Action by $x_1, y_1$ defines linear maps $X_i: V_i\to V_{i+1}$ and $Y_i: V_{i+1}\to V_i$. The space  $\bold{e} E$ can be identified with the subspace of $G$-invariants in $\ee_0 E$, and we have
\begin{equation*}
n\bold{e}=(1+s_{12}+\dots+s_{1n})\ee_0\,.
\end{equation*}
This allows us to define linear maps $v_0: V_\infty\to V_0$ and $w_0: V_0\to V_{\infty}$ by 
\begin{equation}
\label{vw}
v_0(\bold{e} u) = \ee_0 (1+s_{12}+\dots +s_{1n}) u\,,\qquad  w_0(\ee_0 u)= mc_{00}\bold{e} u\,.
\end{equation}
Up to rescaling, these are the canonical inclusion and projection maps. 

We have, for any $\ee_{-i} u\in V_i$,
\begin{multline*}
(Y_iX_i-X_{i-1}Y_{i-1})\ee_{-i} u=[y_1,x_1]\ee_{-i} u=\ee_{-i}\Big(t+c_{00}\sum_{j\ne 1}\sum_{l=0}^{m-1}{s}_{1j}{s}_1^{l}{s}_j^{-l}
+\sum_{l=1}^{m-1}c_l{s}_1^l\Big) u\\
=\ee_{-i} \Big(t+c_{00}\sum_{j\ne 1}\sum_{l=0}^{m-1}{s}_{1j}{\mu}^{-il}
+\sum_{l=1}^{m-1}c_l{\mu}^{-il}\Big) u =\ee_{-i} \Big(t+mc_{00}\delta_{i,0}\sum_{j\ne 1}s_{1j}+\sum_{l=1}^{m-1}c_l{\mu}^{-il}\Big) u\,.
\end{multline*}
It follows that the maps $X_i, Y_i, v_0, w_0$ satisfy relations \eqref{cm0} with 
\begin{equation}\label{para}
\lambda_0=t-mc_{00}+\sum_{l=1}^{m-1} c_l\,,\quad \lambda_i=t+\sum_{l=1}^{m-1}c_l{\mu}^{-il}\ (i\ne 0)\,.
\end{equation}
We also have $w_0v_0=\lambda_\infty\id_{V_\infty}$ with $\lambda_\infty=mnc_{00}$. Thus, to every $H_{t,c}$-module we have associated a representation of $(\Pi^{\llambda})^{\mathrm{opp}}$.

Now let us consider the case $t=0$, denoting $H_{c}=H_{0,c}$. Let $ U_{c}=\bold{e} H_{c} \bold{e} \,$ be 
the corresponding spherical subalgebra of $ H_{c} $. According to \cite[Theorems 1.5, 1.6]{EG}, the algebra $U_{c}$ is finitely-generated commutative, without zero divisors; we can think of it being the ring of regular functions on an affine algebraic variety $\Spec (U_{c})$, with $\bold{e}$ playing the role of the identity. The algebra $U_{c}$ admits a non-commutative deformation $U_{t,c}=\bold{e}H_{t,c}\bold{e}$, which furnishes $U_{c}$ with a Poisson bracket. 

\begin{theorem}[\cite{EG}, Theorem 1.13] For generic $c$, the variety $\Spec (U_{c})$ is isomorphic to the Calogero--Moser space $\CM_{n,\lambda}$ as an algebraic Poisson variety, with the parameters $\lambda_i$ are related to $c$ by \eqref{para} with $t=0$, that is,
\begin{equation}\label{para0}
\lambda_0=-mc_{00}+\sum_{l=1}^{m-1} c_l\,,\quad \lambda_i=\sum_{l=1}^{m-1}c_l{\mu}^{-il}\ (i\ne 0)\,.
\end{equation}
\end{theorem}

\begin{remark} According to \cite[Theorem 7.4]{M}, the genericity assumption on the parameters can be removed, i.e., the varieties $\mathrm{Max}\ \Spec (U_{c})$ and $\CM_{n,\lambda}$ are isomorphic for any $c$.
\end{remark}

Let us recall the construction of the isomorphism between $\Spec (U_{c})$ and $\CM_{n,\lambda}$, following \cite{EG}. First, by Theorems 1.7, 3.1, 16.1 of \cite{EG}, $\Spec (U_{c})$ can be identified with the moduli space of finite-dimensional $H_{c}$-modules. Moreover, for generic $c$ all such modules are irreducible and isomorphic, as $G$-modules, to the regular representation $\C G$. Let $E$ be such a module. As explained above, one can associate to $E$ a representation of the deformed preprojective algebra $(\Pi^{\llambda})^{\mathrm{opp}}$. Since $E\cong \C G$, we easily find that $\dim V_\infty=1$, $\dim V_i=n$. Thus, this represents a point of $\CM_{n,\lambda}$, so we get a map $\theta: \Spec (U_{c})\to \CM_{n,\lambda}$, and as shown in \cite{EG}, $\theta$ is an isomorphism that respects the Poisson brackets.    \qed

According to the above theorem, the induced algebra map $\theta^*: \C  [\CM_{n,\lambda}]\to U_{c}$ is an isomorphism of the coordinate rings. Let us calculate the image of $\tr\,Y^{jm}$ under this isomorphism. According to \eqref{hk} we have 
\begin{equation*}
\tr\,(w_0Y^{jm}v_0)=-\frac{\sum_{i\in I}\lambda_i}{m}\,\tr\,Y^{jm}=c_{00}\,\tr\,Y^{jm}\,,
\end{equation*} 
and we just need to calculate how the map $w_0Y^{jm}v_0$ acts on the one-dimensional space $V_\infty=\bold{e} E$. Applying the definitions \eqref{vw}, we get that 
\begin{equation*} 
w_0Y^{jm}v_0(\bold{e} u)=mc_{00} \bold{e} y_1^{jm}(1+s_{12}+\dots +s_{1n}) u =mc_{00} \bold{e}(y_1^{jm}+\dots +y_n^{jm})(\bold{e}u)\,.
\end{equation*}
Therefore,
\begin{equation*}
\theta^*(\tr \,Y^{jm})=m\,\bold{e}(y_1^{jm}+\dots +y_n^{jm})\bold{e}\in U_{c}\,,\qquad j\in\N\,.
\end{equation*}
The functions $\frac1m \tr\,Y^{jm}$, $j=1, \dots, n$ are in involution on $\CM_{n,\lambda}$ and so they define a completely integrable system. Similarly, the functions $\bold{e}(y_1^{jm}+\dots +y_n^{jm})\bold{e}$ Poisson-commute on the variety $\Spec (U_{c})$ and define a completely integrable system on it. By the above results, the Poisson map $\theta$ establish an isomorphism between these two integrable systems. It remains to explain how this is related to the Calogero--Moser system that is usually constructed with the help of Dunkl operators. This is explained in Section~\ref{secRelDO}.

\subsection{Relation to Dunkl operators}
\label{secRelDO}

A version of the Calogero--Moser system can be defined for any finite reflection group $G$, by using the Dunkl operators \cite{D, DO}; this idea goes back to G.~Heckman \cite{He}. Below we outline that construction in the classical case for $G=\Z_m\wr S_n$, following \cite{Et, EM}. 

For $\h=\C^n$, let $(x_1, \dots, x_n, p_1, \dots, p_n)$ denote the standard canonical coordinates on the symplectic space $\h\times\h^*\simeq T^*(\h)$.
We have an obvious action of $G$ on $\hreg\times\mathfrak h^*$ and on the ring of functions $\C[\hreg\times\mathfrak h^*]$.

Having parameter $c$ as above, let $\kappa_0, \dots , \kappa_{m-1}\in \C$ be chosen in such a way that 
\begin{equation}\label{parac}
\sum_{r=1}^{m-1}c_rs^r=\sum_{l=0}^{m-1}m(\kappa_{l}-\kappa_{l-1})\epsilon_l\,,
\end{equation}
which fixes all $\kappa_i$ uniquely up to a simultaneous shift. The classical Dunkl operators for $G$ are the following elements of $\C[\hreg\times\h^*]*G$:
\begin{equation}\label{didef}
 D_i=p_i-c_{00}\sum_{j\ne i}\sum_{r=0}^{m-1}\frac1{x_i-\mu^rx_j}s_i^r s_{ij} s_i^{-r}
-\frac{1}{x_i}\sum_{l=0}^{m-1}m\kappa_l(\epsilon_{l})_{i}\,,
\end{equation}
where $(\epsilon_{l})_{i}$ denotes the idempotent $\epsilon_l\in\C\Gamma$ placed in the $i$-th component in $\C\Gamma^n$,
\begin{equation*}
(\epsilon_{l})_{i}=\frac{1}{m}\sum_{r=0}^{m-1}\mu^{-lr} s_i^r\,.
\end{equation*}

\begin{prop}[\cite{EM}, Proposition~3.2 and Corollary~3.13] (1) The assignment $x_i\mapsto x_i$, $y_i\mapsto D_i$ and $g\mapsto g$ for $g\in G$ extends to an injective algebra map $\Theta:\, {H}_{c}\rightarrow \C[\hreg\times\h^*]*G$, called the Dunkl embedding. 

(2) When restricted to the spherical subalgebra, the map $\Theta$ defines an embedding $U_{c} \hookrightarrow \bold{e}(\C[\hreg\times\h^*]*G)\bold{e}\simeq \C[\hreg\times\h^*]^G$. Under this embedding, the Poisson structure on $U_{c}$ coming from its non-commutative deformation, $U_{t,c}$, coincides with the restriction of the natural Poisson structure on $\C[\hreg\times\h^*]^G$. 
\end{prop}

Now define a linear map $\bold{m}:\,\C[\hreg\times\h^*]*G\rightarrow \C[\hreg\times\h^*]$ by the formula $\sum_{g\in G} B_g\cdot g\mapsto  \sum_{g\in G} B_g$. Restricting $\bold{m}\circ \Theta$ onto the spherical subalgebra gives the embedding $U_{c} \hookrightarrow \C[\hreg\times\h^*]^G$, as in the proposition above. The following result is a special case of the general construction of the classical Calogero--Moser system for an arbitrary complex reflection group, see~\cite[Section 6.3]{Et} and~\cite[Section 2.10]{EM}. 

\begin{prop}\label{cmint} The functions $L_j(x,p)=\bold{m}\circ \Theta ( \bold{e} (y_1^{jm}+\dots +y_n^{jm})\bold{e})$, $j=1,\dots, n$, are independent $G$-invariant and Poisson-commuting functions on $\hreg\times\h^*=T^*\hreg$. Hence, they define a completely integrable system on $T^*\hreg$.
\end{prop}

\begin{remark} According to \cite[Lemma 2.2]{EFMV}, $\Theta ( \bold{e} (y_1^{jm}+\dots +y_n^{jm})\bold{e})$ is already in $\C[\hreg\times\h^*]$, i.e., it does not contain nontrivial elements of $G$. Thus, application of $\bold{m}$ in the construction of $L_j$ is not necessary.  
\end{remark}

\begin{remark} \label{remKappa}
It is customary to assume that $\sum_{l=0}^{m-1}\kappa_l=0$ in the definition of the Dunkl operators $D_i$. In this case the Hamiltonians $L_j(x,p)$ have the form $\sum_{i=1}^n p_i^{mj}$ plus the terms of degree $\le mj-2$ in momenta (cf. Example~\ref{ExCMBn}).
\end{remark}

Earlier we have found  that $\bold{e}L_j(x,p)\bold{e}$ equals $\frac1m \theta^*(\tr\, Y^{jm})$. Now we will make this relation more explicit, following the approach of \cite[Section 9.6]{EM}. Recall that $\CM_{n,\lambda}$ has an open dense subset $\CM_{n,\lambda}'$ with local coordinates $\{\mu_i, \nu_i\}$, see \eqref{slice0} and \eqref{slice1}--\eqref{slice3}. We are going to construct similar coordinates on a subset of $\Spec (U_{c})$. Take ${\mu}=({\mu}_1,\ldots,{\mu}_n)\in\mathfrak h^*$, ${\nu}=({\nu}_1,\ldots,{\nu}_n)\in\hreg$, and consider $\mathbb C_{{\nu},{\mu}}=\mathbb C\cdot 1_{{\nu},{\mu}}$ to be a one-dimensional representation of the algebra $\C[\hreg\times\h^*]$ defined by
\begin{equation*}
x_i\cdot1_{{\nu},{\mu}}={\nu}_i\cdot1_{{\nu},{\mu}}\,,\qquad p_i\cdot1_{{\nu},{\mu}}={\mu}_i\cdot1_{{\nu},{\mu}}\,.
\end{equation*}
Viewing $\C[\hreg\times\h^*]$ as a subalgebra of $\C[\hreg\times\h^*]*G$, we take $E$ to be the representation of $\C[\hreg\times\h^*]*G$ induced from $\C_{\nu, \mu}$. Composing this with the Dunkl embedding $\Theta$ gives us a representation, denoted $E_{\nu, \mu}$, of $H_{c}$. Clearly, $E_{\nu, \mu}\simeq \C G$ as a $G$-module, so it is spanned by the elements $g \cdot 1_{\nu, \mu}$, $g\in G$.
It was explained above how to associate to $E_{\nu, \mu}$ a point of $\CM_{n,\lambda}$. We have $V_i=\ee_{-i}E_{\nu, \mu}$, and we can choose the elements $\ee_{-i}\,s_{1j} \cdot1_{\nu, \mu}$, $j=1, \dots, n$, as a basis of $V_i$, identifying it with $\C^n$. Here we use the convention that $s_{11}=1$. Also, we have $V_{\infty}=\bold{e}E_{\nu, \mu}=\C\bold{e} \cdot 1_{\nu, \mu}$, which we identify with $\C$. With these identifications, the linear maps $X_i, Y_i, v_0, w_0$ become matrices and column/row vectors, respectively. 

We have $x_1(\ee_{-i}\,s_{1j} \cdot1_{\nu, \mu})=\ee_{-i-1}\,s_{1j} x_j \cdot1_{\nu, \mu}=\nu_j\ee_{-i-1}\,s_{1j} \cdot1_{\nu, \mu}$. Therefore, the matrices $X_i$ look exactly like in \eqref{slice0}. Also, from the definitions of $v_0, w_0$ and with the help of \eqref{para0} it is easy to check that \eqref{slice1} holds. It remains to calculate $Y_i$. For example, to calculate the matrix entries of $Y_0$ we consider
$y_1(\ee_{-1}\,s_{1i} \cdot1_{\nu, \mu})=\ee_{0}\,s_{1i}D_i\cdot1_{\nu, \mu}$.
Using that 
\begin{equation*}
D_i=p_i+c_{00}\sum_{j\ne i}\sum_{r=0}^{m-1}s_i^r s_{ij} s_i^{-r}(x_i-\mu^rx_j)^{-1}
-\sum_{l=0}^{m-1}m\kappa_l(\epsilon_{l+1})_{i}x_i^{-1}\,,
\end{equation*}
we find that
\begin{equation*}
y_1(\ee_{-1}\,s_{1i} \cdot1_{\nu, \mu})=(\mu_i-\frac{m\kappa_{-1}}{\nu_i})\ee_0\, s_{1i}\cdot1_{\nu, \mu} +c_{00}\sum_{j\ne i}\sum_{r=0}^{m-1}\frac{1}{\nu_i-\mu^r\nu_j}\ee_0\,s_{1j}\cdot1_{\nu, \mu}\,.
\end{equation*}   
This agrees with the formulas \eqref{slice2}--\eqref{slice3} for $Y_0: V_1\to V_{0}$, provided that $a+\lambda_0=-m\kappa_{m-1}$, that is,
\begin{align}\label{match}
  \kappa_{m-1}=\sum_{r=1}^{m-1}\frac{m-r}{m^2}\lambda_r\,.
\end{align}
This together with~\eqref{parac} fixes uniquely the parameters $\kappa_l$. Note that once $Y_0$ is known, all other matrices $Y_i$ are uniquely determined from the relations \eqref{mr1}--\eqref{mr2}. Therefore, they will also agree with \eqref{slice2}--\eqref{slice3}. Hence, the representations $E_{\nu, \mu}$ with $(\nu,\mu)\in\hreg\times\mathfrak h^*$ form a dense subset of $\Spec (U_{c})$, with the coordinates $\{\nu_i, \mu_i\}$ matching exactly those on $\CM_{n,\lambda}$ under the isomorphism $\theta: \Spec (U_{c})\to \CM_{n,\lambda}$. We may think of $\Spec (U_{c})$ as the moduli space of one-dimensional representations of $U_{c}$: given $E=E_{\nu, \mu}$ as above, the corresponding one-dimensional representation is $\bold{e} E_{\nu, \mu}=\C\bold{e}\cdot 1_{\nu, \mu}$. A function $\bold{e}L_j(x,p)\bold{e}=\frac1m \theta^*(\tr\, Y^{jm})$ acts on $\bold{e}\cdot 1_{\nu, \mu}$ by left multiplication, resulting in $\bold{e}L_j(\nu,\mu)\cdot 1_{\nu,\mu}$. This implies the following results.

\begin{prop}\label{cmcn} Given arbitrary $c$, define $\lambda_i$ and $\kappa_l$ by \eqref{para0}, \eqref{parac} and \eqref{match}. Let the matrices $Y_i$ be as in \eqref{slice2}--\eqref{slice3}, and $Y=Y(\nu, \mu)$ be the $mn\times mn$ matrix $Y=\sum_i Y_i\in \End(\oplus_i V_i)$. Then for any $j\in\N$ we have $\frac1m \tr ( Y^{jm})=L_j(\nu, \mu)$, where $L_j=L_j(x,p)$ are the Hamiltonians of the Calogero--Moser problem as in Prop.~\ref{cmint}.
\end{prop}  

\begin{cor} The matrix $Y=\sum_{i} Y_i$ is a Lax matrix for the Calogero--Moser system for $G=\Z_m\wr S_n$, in the sense that it deforms isospectrally under any of the commuting flows defined by the Hamiltonians $L_j(\nu, \mu)$, and the traces of powers of $Y$ produce the full set of first integrals in involution. The symplectic quotient $\CM_{n,\lambda}$ is a completed phase space for the Calogero--Moser system, with the flows extending to $\CM_{n,\lambda}$ and obtained by the symplectic reduction from the linear flows \eqref{hd}.
\end{cor}

\begin{remark}
 Using the formulas $m(\kappa_{l}-\kappa_{l-1})=\sum_{r=1}^{m-1}c_r\mu^{rl}=\lambda_{-l}-\delta_{l0}|\lambda|$ one can check that the relation~\eqref{match} leads to $\sum_{l=0}^{m-1}\kappa_l=0$, which agrees with Remarks~\ref{remKappa} and~\ref{shi}.
\end{remark}

\begin{example} \label{ExCMBn} \normalfont
 Consider the case $m=2$. In this case $c=(c_{00},c_1)$, $\lambda_0=-2c_{00}+c_1$, $\lambda_1=-c_1$. This gives $a=2c_{00}-c_1/2$ and $\kappa_0=-\kappa_{1}=c_1/4$. We will use $x,p$ instead of $\nu,\mu$. The matrix $Y$ then has the form
$\begin{pmatrix}0&Y_1\\ Y_0 &  0 \end{pmatrix}$ where
\begin{align*}
&(Y_0)_{ii}=p_i+\frac{c_1}{2x_i}\,,\qquad (Y_0)_{ij}=\frac{-2c_{00}x_j}{x_i^2-x_j^2}\,,\\
&(Y_1)_{ii}=p_i-\frac{c_1}{2x_i}\,,\qquad (Y_1)_{ij}=\frac{-2c_{00}x_i}{x_i^2-x_j^2}\,.
\end{align*}
Up to a change of basis, this is equivalent to the well-known Lax matrix for the rational Calogero--Moser system in type~$B_n$~\cite{OP1}, see also~\cite[eq.~(4.48)]{BCS}. We have
\begin{align*}
\frac12 \tr\, Y^2= \sum_{i=1}^n\left(p_i^2-\frac{c_1^2}{4x_i^2}\right)-\sum_{i< j}\frac{4c_{00}^2(x_i^2+x_j^2)}{(x_i^2-x_j^2)^2}\,.
\end{align*}
Note that setting $c_1=0$ leads to the $D_n$ case. 
\end{example}

\begin{remark}
 Note that the link with the Calogero--Moser system was established only for a special choice of the dimension vector $\aalpha=(1,n\delta)$. It would be interesting to understand the situation for general $\aalpha=(1,\alpha)\in\Delta_+(Q)$. We expect that every imaginary positive root $(1,\alpha)$ can be reduced to $(1,n\delta)$ by a sequence of simple reflections. In such a case the reflection functors~\cite{CB2} induce an isomorphism of the symplectic varieties $\CM_{\alpha,\lambda}$ and $\CM_{n,\lambda}=\CM_{n\delta,\lambda'}$ for generic $\lambda$. Moreover, one can also check that the reflection functors preserve the Poisson-commutative subalgebra generated by $\tr\,Y^{mj}$, $j\ge1$.
\end{remark}

\subsection{Spin particle dynamics}
\label{secFullHier}

For the full hierarchy~\eqref{kp}, the solutions are parameterised by the varieties \eqref{cms} for the quiver $Q$ as in \eqref{q}. As in the previous section, let us consider the special choice of the dimension vector $\aalpha=(1, n\delta)$ where $\delta=(1, \dots , 1)$. Let us again abbreviate the varieties $\CM_{n\delta, \lambda}(Q)$ to $\CM_{n,\lambda}$ for simplicity. From the formula~\eqref{dimq}, we find that $\dim \CM_{n,\lambda}=2mn$. Points of $\CM_{n,\lambda}$ are represented by linear maps
\begin{align*}
 &X_{i}\colon V_i\to V_{i+1}\,, &&Y_{i}\colon V_{i+1}\to V_i\,, &&v_{i}\colon \C\to V_i\,, &&w_{i}\colon V_i\to \C\,,    
\end{align*}
satisfying the relations~\eqref{cm1}:
\begin{align}
\label{cmg}
 &Y_{i}X_{i}-X_{i-1}Y_{i-1} - v_{i}w_{i}=\lambda_i\id_{V_i}  &&(i\in I)\,.
\end{align}
Similarly to the previous case, we look at the subset $\CM'_{n,\lambda}\subset \CM_{n,\lambda}$ where $X_{m-1}\dots X_1X_0$ has distinct non-zero eigenvalues $\nu_1^m,\ldots,\nu_n^m$. Then by choosing suitable bases in each of $V_i$ we will have
\begin{equation}
\label{sliceg}
X_i=\mathrm{diag}({\nu}_1, \dots, {\nu}_n)\,,
\end{equation}
implying  that
\begin{align}
\label{mrg1}
(Y_i)_{jj}\nu_j-\nu_j(Y_{i-1})_{jj}-(v_iw_i)_{jj}&=\lambda_i\,,\\
\label{mrg2}
(Y_i)_{jk}\nu_k-\nu_j(Y_{i-1})_{jk}-(v_iw_i)_{jk}&=0 &&(j\ne k)\,.
\end{align}
Summing~\eqref{mrg1} over $i$ we obtain that for any $j$, 
\begin{equation*}
\sum_{i\in I}(v_iw_i)_{jj}=-|\lambda|\,.
\end{equation*} 
The matrices $Y_i$ are then found to be as follows  
\begin{align}
\label{slice2g}
&(Y_i)_{jj}=\mu_j+\frac{1}{\nu_j}\sum_{r=0}^{i}\big(\lambda_r+(v_rw_r)_{jj}\big)+\frac{1}{\nu_j}\sum_{r=0}^{m-1}\frac{r-m}m\big(\lambda_r+(v_rw_r)_{jj}\big)\,,\\
 \label{slice3g}
&(Y_i)_{jk}=-\sum_{l=0}^{m-1}\,\frac{(v_{i-l}w_{i-l})_{jk}\nu_j^{l}\nu_k^{m-l-1}}{\nu_j^m-\nu_k^m}\qquad\qquad (j\ne k)\,,
\end{align}
where $i=0, 1, \dots, m-1$ (the index $i-l$ in the second formula is taken modulo $m$). 

\begin{remark}
 The second sum in the right hand side of~\eqref{slice2g} is added to have the relation $\sum_{i\in I}(Y_i)_{jj}=m\mu_j$. This choice is important for Prop.~\ref{sfs} below.
\end{remark}

It seems convenient to organise the column and row vectors $v_i, w_i$ into matrices $v$ and $w$ of size $n\times m$ and $m\times n$, respectively, and then divide them into rows and columns, each of length $m$. This way we get $n$ rows $\varphi_j$ and $n$ columns $\psi_j$, with $(v_i)_j=(\varphi_j)_i$ and $(w_i)_j=(\psi_j)_i$. The formulas for $Y_i$ then become as follows:
\begin{align}
\label{slice2gg}
&(Y_i)_{jj}=\mu_j+\frac{1}{\nu_j}\sum_{r=0}^{i}\big(\lambda_r+(\psi_j\varphi_j)_{rr}\big)+\frac{1}{\nu_j}\sum_{r=0}^{m-1}\frac{r-m}m \big(\lambda_r+(\psi_j\varphi_j)_{rr}\big) \,,\\
 \label{slice3gg}
&(Y_i)_{jk}=-\sum_{l=0}^{m-1}\,\frac{(\psi_k\varphi_j)_{i-l,i-l}\,\nu_j^{l}\nu_k^{m-l-1}}{\nu_j^m-\nu_k^m}\qquad\qquad (j\ne k)\,,
\end{align}
with the constraint $\varphi_j\psi_j = -|\lambda|$ for all $j$. Note that acting on each $V_i$ by a diagonal matrix with diagonal entries $d_1, \dots, d_n$ in accordance with \eqref{act}, leaves $\nu_j, \mu_j$ unchanged and replaces $\varphi_j$, $\psi_j$ by $d_j\varphi_j$ and $d_j^{-1}\psi_j$. By analogy with \cite{GH} (cf.~\cite{W2}), the quantities $\varphi_j, \psi_j$ can be thought of as spin variables attached to the $j$-th particle of an $n$-particle system of Calogero--Moser type. With each pair $(\varphi_j, \psi_j)$ we can therefore associate a point $q_j$ of the symplectic variety $Q_m$ (cf.~\cite{W2}):
\begin{equation*}
Q_m=\big\{(\varphi, \psi)\in(\C^m)^*\oplus \C^m \mid \varphi \psi=-|\lambda|\big\}/\C^\times\,.
\end{equation*}  
The above formulas then give a local parameterisation of $\CM_{n,\lambda}'$ by $\{\nu_j, \mu_j, q_j\}$.

By an elementary residue calculation, 
\begin{equation*}
\frac{\nu_j^{l}\nu_k^{m-l-1}}{\nu_j^m-\nu_k^m}=\frac{1}{m}\sum_{r=0}^{m-1}\frac{\mu^{-lr}}{\mu^{r}\nu_j-\nu_k}\,. 
\end{equation*}
As a result, the formula \eqref{slice3gg} can be replaced with
\begin{equation}\label{slice4gg}
(Y_i)_{jk}=\sum_{r\in I} \,\frac{\mu^{-ir}\left(\varphi_j \, E^r \,\psi_k\,\right)}{m(\nu_k-\mu^{r}\nu_j)}\,,\qquad E^r=\mathrm{diag}(1, \mu^{r}, \mu^{2r}, \dots, \mu^{(m-1)r})\,.
\end{equation} 

Recall that the symplectic structure \eqref{symp1} induces a symplectic structure on $\CM_{n,\lambda}$.

\begin{prop} \label{sfs}
 In terms of $\nu_j, \mu_j, q_j=(\varphi_j, \psi_j)$ the symplectic structure on $\CM_{n,\lambda}$ has the form
\begin{equation*}
 \omega_{\CM_{n,\lambda}}=\sum_{j=1}^n (m\,\mathrm{d}\mu_j \wedge \mathrm{d}\nu_j-\mathrm{d}\varphi_j \wedge \mathrm{d}\psi_j)\,. 
\end{equation*}
Here, by abuse of notation, we do not make distinction between the form $\mathrm{d}\varphi_j \wedge \mathrm{d}\psi_j=\sum_{r=0}^{m-1}\mathrm{d}\varphi_{j,r} \wedge \mathrm{d}\psi_{j,r}$ on $T^*\C^m$ and the corresponding symplectic form on~$Q_m$.
\end{prop}

\noindent{\bf Proof.} The symplectic form on $\Rep\big(\C\overline{Q}^{\mathrm{opp}},\aalpha\big)$ in this case is $$\omega=\sum\limits_{i\in I}\Big(\tr(\mathrm{d}Y_i \wedge \mathrm{d}X_i)+\mathrm{d}w_i \wedge \mathrm{d}v_i\Big)\,.$$ Let view the formulas~\eqref{sliceg}, \eqref{slice2gg}, \eqref{slice3gg} as describing a map $$\phi\colon T^*\hreg\times(T^*\C^m)^n\to\Rep\big(\C\overline{Q}^{\mathrm{opp}},\aalpha\big)\,.$$ Similarly to the proof of Prop.~\ref{sf} we need to calculate $\phi^*\omega$. In the same way, only the diagonal entries of $Y_i$ will contribute to $\phi^*\omega$. We have
\begin{align*}
 \mathrm{d}(Y_i)_{jj}=\mathrm{d}\mu_j+\frac1{\nu_j}\Big(\sum_{r=0}^i\mathrm{d}\big(\psi_j\varphi_j\big)_{rr}+\sum_{r=0}^{m-1}\frac{r-m}m\mathrm{d}\big(\psi_j\varphi_j\big)_{rr}\Big)+\frac{\partial(Y_i)_{jj}}{\partial\nu_j}\mathrm{d}\nu_j\,.
\end{align*}
The last term will vanish after taking the wedge product with $\mathrm{d}(X_i)_{jj}=\mathrm{d}\nu_j$. A simple calculation then shows that $\sum_{i=0}^{m-1}\mathrm{d}(Y_i)_{jj} \wedge \mathrm{d}(X_i)_{jj}=m\,\mathrm{d}\mu_j \wedge \mathrm{d}\nu_j$. \qed
 
The Hamiltonians $H_k$ \eqref{ham} after substituting \eqref{slice2gg}--\eqref{slice3gg} produce Poisson commuting functions on the phase space $\hreg \times \mathfrak h^*\times (Q_m)^n$ equipped with the symplectic form $\omega_{\CM_{n,\lambda}}$.
They define a system which can be interpreted as a spin version of the $n$-particle Calogero--Moser system for $G=\Z_m\wr S_n$, with the spin variables $(\varphi_j, \psi_j)\in Q_m$. It is a spacial case of a more general system with spin variables, see Subsection~\ref{secDarCoorMat}.

Note that the choice of coordinates on $\CM'_{n,\lambda}$ is only unique up to the action of $G$, coming from the embedding~\eqref{embed}. Explicitly in coordinates this action permutes $\{\nu_j\}$, $\{\mu_j\}$, $\{q_j=(\varphi_j,\psi_j)\}$ simultaneously and replaces  $\nu_j$, $\mu_j$, $\varphi_j$, $\psi_j$ by $\nu^{r_j}\nu_j$, $\mu^{-r_j}\mu_j$, $\varphi_j E^{r_j}$, $E^{-r_j}\psi_j$, respectively, where $E^r$ is the same as in~\eqref{slice4gg}. It follows that $\CM'_{n,\lambda}\simeq\hreg \times \mathfrak h^*\times (Q_m)^n/G$. We obtain the following result.

\begin{prop} Let $H_k=H_k(\nu, \mu, q)$ be the functions of $(\nu, \mu, q)\in \hreg \times \mathfrak h^*\times (Q_m)^n$, obtained by substituting \eqref{slice2gg}--\eqref{slice3gg} into \eqref{ham}. Then $H_k$ with $k=1, \dots, mn$ are functionally independent, Poisson commuting functions and hence define a completely integrable system. The Hamiltonian flows defined by $H_k$ extend from $\hreg \times \mathfrak h^*\times (Q_m)^n/G$ to produce complete flows on the Calogero--Moser space $\CM_{n,\lambda}=\CM_{n\delta, \lambda}$ for the quiver \eqref{q}.  
\end{prop}

The number of integrals, $mn$, is clearly the half of the dimension of the phase space $\hreg \times \mathfrak h^*\times (Q_m)^n$, so the only thing to prove is that $H_k$ with $k=1, \dots, mn$ are independent. This is shown in a more general case in Section~\ref{secMultGen} (Proposition~\ref{liouv}). \qed 

\begin{example} \label{m=2}
\normalfont
Consider the case $m=2$. Let us write down the first two Hamiltonians in the coordinates $\nu_j, \mu_j, q_j$. In this case we have two-component spin variables $\varphi_j=(\varphi_{j0}, \varphi_{j1})$, $\psi_j=(\psi_{j0}, \psi_{j1})^\top$. 
Denote
$E=\begin{pmatrix}1 & 0 \\ 0 & -1\end{pmatrix}$, $F_\pm=\begin{pmatrix}0 & \pm1 \\ 1 & 0\end{pmatrix}$.
Then we have
\begin{align*}
H_1&=-\sum_{i=0,1}w_iY_iv_{i+1}=-\sum_{j=1}^n\Big((\varphi_jF_+\psi_j)\mu_j+\frac1{2\nu_j}(\varphi_j F_-\psi_j)(\lambda_0+\varphi_{j0}\psi_{j0})\Big)\\
 &\qquad\qquad\qquad\qquad
-\sum_{j\ne k}\Big(\frac{(\varphi_j F_+\psi_k)(\varphi_k\psi_j)}{2(\nu_j-\nu_k)}
+\frac{(\varphi_j F_-\psi_k)(\varphi_k E\psi_j)}{2(\nu_j+\nu_k)}\Big)\,, \\
H_2&=\frac{|\lambda|}2\tr Y^2=|\lambda|\sum_{j=1}^n\Big(\mu_j^2-\frac1{4\nu_j^2}(\lambda_0+\varphi_{j0}\psi_{j0})^2\Big)\\
 &+|\lambda|\sum_{j\ne k}\Big(\frac{(\varphi_j\psi_k)(\varphi_k\psi_j)}{4(\nu_j-\nu_k)^2}
+\frac{(\varphi_j E\psi_k)(\varphi_k E\psi_j)}{4(\nu_j+\nu_k)^2}
+\frac{\varphi_{j1}\psi_{k1}\varphi_{k0}\psi_{j0}-\varphi_{j0}\psi_{k0}\varphi_{k1}\psi_{j1}}{2(\nu_j^2-\nu_k^2)}\Big)\,,
\end{align*}
where $(\lambda_0+\varphi_{j0}\psi_{j0})=-(\lambda_1+\varphi_{j1}\psi_{j1})$.
\end{example}

In the case $m=2$, when $G$ is of $B_n$ type, some spin versions of Calogero-Moser problem can be constructed in the framework of Lie algebras, see~\cite{LX} and references therein. However they look somewhat different to our system. To the best of our knowledge, for $m>2$ spin generalizations of the Calogero--Moser system for $G=\Z_m\wr S_n$ have not been considered before.

\section{Multicomponent generalisation}
\label{secMultGen}

The KP hierarchy admits a multicomponent (matrix) generalisation~\cite{DJKM}, \cite{Sa}, see also~\cite{Di}. In the same spirit one can generalise the hierarchy~\eqref{kp}. Below we introduce such a hierarchy, and construct its solutions by using quivers.

\subsection{Definition of the hierarchy}
\label{secDefHier}

We will follow closely the notation of~\cite{Di}.

We will work with the ring $\mathop{\mathrm{Mat}}(d,\mathbb C)\otimes\P$ of matrices over the ring $\P$ \eqref{p}. The operations $(\cdot)_+$ and $(\cdot)_-$ on $\P$ induce similar operations on this ring.

Let $E_\alpha$ ($\alpha=1,\ldots,d$) be the diagonal matrices with the diagonal entries $(E_\alpha)_{\beta\beta}=\delta_{\alpha\beta}$. We have $E_\alpha E_\beta=\delta_{\alpha\beta}E_\alpha$ and $E_1+\ldots+E_d=1$.

Now consider the elements of $\mathop{\mathrm{Mat}}(d,\mathbb C)\otimes\P$ of the following form
\begin{align}
 &L=y+\sum_{l=0}^\infty f_ly^{-l}, \label{Lmat} \\
 &R_\alpha=E_\alpha+\sum_{l=1}^\infty r_{l\alpha}y^{-l} &&(\alpha=1,\ldots,d) \label{Ralpha}
\end{align}
where $f_l,r_{l\alpha}\in\mathop{\mathrm{Mat}}(d,\mathbb C)\otimes\big(\mathbb C(x)*\Gamma\big)$.

We impose the following relations on $L$, $R_\alpha$:
\begin{align}
 &[L,R_\alpha]=0\,, &&R_\alpha R_\beta=\delta_{\alpha\beta}R_\alpha\,, &&\sum_{\alpha=1}^d R_\alpha=1\,. \label{LRalpha}
\end{align}
The matrix generalisation of the hierarchy~\eqref{kp} is the following system of equations on $\big\{L,R_\alpha\big\}$:
\begin{align}
 &\frac{\partial L}{\partial t_{k\beta}}=\big[(L^kR_\beta)_+,L\big]\,, \label{tkbetaL} \\
 &\frac{\partial R_\alpha}{\partial t_{k\beta}}=\big[(L^kR_\beta)_+,R_\alpha\big]\,, \label{tkbetaRalpha}
\end{align}
where $t_{k\beta}$ ($\beta=1,\ldots,d$, $k=0,1,\ldots$) are the time variables. Note that it follows from~\eqref{LRalpha} that $\sum_{\beta=1}^d \frac{\partial}{\partial t_{0\beta}}=0$.

The following proposition and its proof is parallel to Prop.~\ref{prop2.1}

\begin{prop} \label{prop7.1}

{\normalfont (1)} The equations~\eqref{tkbetaL}, \eqref{tkbetaRalpha} imply the zero-curvature equations:
\begin{align}
 \Big[\frac{\partial}{\partial t_{k\beta}}-(L^kR_\alpha)_+\;,\;\frac{\partial}{\partial t_{l\beta}}-(L^lR_\beta)_+\Big]=0 &&&(k,l\ge0). \label{zccmat}
\end{align}
{\normalfont (2)} Flows~\eqref{tkbetaL}, \eqref{tkbetaRalpha} pairwise commute and preserve the relations~\eqref{LRalpha}. 
\end{prop}

We will also impose the following constraint analogous to~\eqref{f0}:
\begin{align}
 &(f_{0,\alpha\alpha})_e=0 &&(\alpha=1,\ldots,d), \label{f0mat}
\end{align}
where $f_{0,\alpha\alpha}\in\mathbb C(x)*\Gamma$ are the diagonal entries of the coefficient $f_0$ in~\eqref{Lmat}. These conditions~\eqref{f0mat} are preserved by the flows~\eqref{tkbetaL}, \eqref{tkbetaRalpha}.

For $M\in1+\mathop{\mathrm{Mat}}(d,\mathbb C)\otimes\P_-$, define $L$ and $R_\alpha$ by the following formulas:
\begin{align}
 &L=MyM^{-1}, &&R_\alpha=ME_\alpha M^{-1}. \label{LM}
\end{align}
It is clear that these elements have the form~\eqref{Lmat}, \eqref{Ralpha} and satisfy~\eqref{LRalpha}, \eqref{f0mat}. Note that $L^kR_\alpha=My^kE_\alpha M^{-1}$.

\begin{prop} \label{propMmat} {\normalfont (cf. Prop.~\ref{wave})}. Suppose $M\in1+\mathop{\mathrm{Mat}}(d,\mathbb C)\otimes\P_-$ depends on $t_{k\beta}$ and satisfies the equation
\begin{align}
 \frac{\partial M}{\partial t_{k\beta}}=-\big(My^kE_\beta M^{-1}\big)_-\cdot M\,. \label{Mmat}
\end{align}
Then $L$ and $R_\alpha$ defined by~\eqref{LM} satisfy the equations~\eqref{tkbetaL}, \eqref{tkbetaRalpha}.
\end{prop}

\begin{remark}\label{remLA}
 A more general hierarchy can be constructed by allowing in~\eqref{LRalpha}--\eqref{tkbetaRalpha} $L$ of the form
\begin{align}
 &L=Ay+\sum_{l=0}^\infty f_ly^{-l}, \label{LmatA}
\end{align}
where $A=\diag(a_1,\ldots,a_d)$ with non-zero $a_\alpha\in\mathbb C$ (cf.~\cite{Dickey}). Solutions to this hierarchy can then be constructed as follows. Let $M$ be an element satisfying the assumptions of Prop.~\ref{propMmat}. Then $L_A=MAyM^{-1}$ and $R_\alpha=ME_\alpha M^{-1}$ satisfy the equations~\eqref{LRalpha}--\eqref{tkbetaRalpha} after rescaling the time variables by $t_{k\beta}\to a_\beta^{-k}t_{k\beta}$.
\end{remark}

\subsection{Constructions of solutions}
\label{secConSol}

We are going to construct solutions of the hierarchy~\eqref{LRalpha}--\eqref{tkbetaRalpha}, parametrised by Calogero--Moser spaces associated with the quivers~\eqref{qd}, \eqref{qld}.

Let $Q$ be the quiver~\eqref{qd} and $\CM_{\alpha,\lambda}=\CM_{\alpha,\lambda}(Q)$ be the Calogero--Moser space~\eqref{cms}:
\begin{align}
\label{cmsmat}
 &\CM_{\alpha, \lambda}(Q)=\mathtt{Rep}\big((\Pi^{\llambda})^{\mathrm{opp}}, \aalpha\big)//\mathtt{GL}(\alpha)\,. 
\end{align}
Recall that here $\alpha\in\N_0^m$, $\lambda\in\C^m$ and $\aalpha=(1,\alpha)$, $\llambda=(-\lambda\cdot\alpha,\lambda)$. A point of $\CM_{\alpha,\lambda}$ represented by a collection of maps
\begin{align*}
 &X_i\colon V_i\to V_{i+1}\,, & &Y_i\colon V_{i+1}\to V_i\,, \\
 &v_{i,\alpha}\colon \mathbb C\to V_i\,, & &w_{i,\alpha}\colon V_i\to\mathbb C & &&(\alpha=1,\ldots,d)
\end{align*}
corresponding to arrows $a_i$, $a_i^*$, $b_{i,\alpha}$ and $b_{i,\alpha}^*$, respectively. These maps should satisfy the relations (cf.~\eqref{cm1}--\eqref{cm2})
\begin{align}
 Y_iX_i-X_{i-1}Y_{i-1}-\sum_{\alpha=1}^d v_{i,\alpha}w_{i,\alpha}&=\lambda_i 1_{V_i}\,, \label{XYvwkmatalpha} \\
 \sum_{i\in I}\sum_{\alpha=1}^d w_{i,\alpha}v_{i,\alpha}&=\lambda_\infty\,. \label{wkvkmat}
\end{align}
(As before, \eqref{wkvkmat} follows from~\eqref{XYvwkmatalpha}.)

Let us define maps
\begin{align*}
 &v_i\colon\mathbb C^d\to V_i\,, & &w_i\colon V_i\to\mathbb C^d
\end{align*}
by $v_i=(v_{i,1},\ldots,v_{i,d})$, $w_i=(w_{i,1},\ldots,w_{i,d})$. Then the relations~\eqref{XYvwkmatalpha} take the form
\begin{align}
 &Y_iX_i-X_{i-1}Y_{i-1}-v_iw_i=\lambda_i 1_{V_i}\,, \label{XYvwkmat}
\end{align}
looking similar to~\eqref{cm1}.

Introducing $V=V_0\oplus\ldots\oplus V_{m-1}$ and $X,Y\in\End(V)$ as in~\eqref{xy} the relation~\eqref{mom} still holds true.

We will be assuming that $\lambda$ is regular in the sense of~\eqref{lambdareg}. Then using Prop.~\ref{dim} we find, similarly to~\eqref{dimq}, that
\begin{align}
 \dim \CM_{\alpha,\lambda}=2\sum_{i\in I}(\alpha_i\alpha_{i+1}+d\alpha_i-\alpha_i^2)=2d\sum_{i\in I}\alpha_i-\sum_{i\in I}(\alpha_i-\alpha_{i+1})^2. \label{dimqmat}
\end{align}

Recall that the space $\CM_{\alpha,\lambda}$ is a symplectic quotient (see~\eqref{sympq}). The symplectic structure on it is induced by the symplectic structure on $\mathtt{Rep}\big((\C Q)^{\mathrm{opp}}, \aalpha\big)$ which has the form
\begin{align}
\label{symp1mat}
\omega=\sum_{i\in I}\big(\tr(\mathrm{d}Y_i \wedge \mathrm{d}X_{i})+\tr(\mathrm{d}w_i \wedge \mathrm{d}v_i)\big).
\end{align}
The Poisson bracket corresponding to~\eqref{symp1mat} are
\begin{align}
\label{poissonmat}
 &\big\{(X_i)_{ab}, (Y_j)_{cd}\big\}=\delta_{ij}\delta_{ad}\delta_{bc}\,, &&\big\{(v_{i,\alpha})_{a}, (w_{j,\beta})_{b}\big\}=\delta_{ij}\delta_{ab}\delta_{\alpha\beta}\,.
\end{align}

Following the idea of~\cite{GH}, consider the following functions on $\CM_{\alpha,\lambda}$:
\begin{align*}
 &I_k(A):=\sum_{i\in I}\tr\big(Aw_i Y^k v_{i+k}\big), &&A\in\mathop{\mathrm{Mat}}(d,\mathbb C).
\end{align*}
The functions $I_k(A)$ with $k=0,1,\ldots$ and $A\in\mathop{\mathrm{Mat}}(d,\mathbb C)$ form a vector space, closed under the Poisson bracket:
\begin{align*}
 &\big\{I_k(A),I_l(B)\big\}=I_{k+l}\big([A,B]\big).
\end{align*}
Choosing $A=-E_\beta$ ($\beta=1,\ldots,d$), we get Poisson-commuting functions $H_{k,\beta}:=-I_k(E_\beta)$. The Hamiltonian flow on $\CM_{\alpha,\lambda}$ looks as follows:
\begin{align}
&\dot X=-\sum_{l=0}^{m-1}\sum_{a=0}^{k-1}Y^av_lE_\beta w_{l-k} Y^{k-a-1}\,, &
&\dot Y=0\,, \label{tFlowXYmat} \\
&\dot w_i=E_\beta w_{i-k}Y^k\,, &
&\dot v_j=-Y^kv_{j+k}E_\beta\,. \label{tFlowVWmat}
\end{align}

Note that these equations can be integrated explicitly similarly to~\eqref{vwsol},~\eqref{Xsol}.

To construct solutions of the hierarchy~\eqref{tkbetaL},\eqref{tkbetaRalpha} let us introduce an element $M\in\mathop{\mathrm{Mat}}(d,\mathbb C)\otimes\P$ by the same formula~\eqref{m}:
\begin{align}
\label{mmat}
M=1-\sum_{i, j\in I}\epsilon_i w_i (X-x\id_V)^{-1}(Y-y\id_V)^{-1} v_j\epsilon_j\,.
\end{align}

\begin{prop} \label{prop6.4}
 Let $L=MyM^{-1}$ and $R_\alpha=ME_\alpha M^{-1}$ with $M$ given by~\eqref{mmat}. If $X,Y,v_i,w_i$ satisfy the equations~\eqref{tFlowXYmat}, \eqref{tFlowVWmat} then $L$ and $R_\alpha$ ($\alpha=1,\ldots,d$) satisfy the equations $\dot L=[(L^kR_\beta)_+,L]$, $\dot R_\alpha=[(L^kR_\beta)_+,R_\alpha]$.
\end{prop}

\noindent{\bf Proof.} Proof essentially repeats the arguments from Section~\ref{secRatSol}.

First, since the relations for $X,Y,v_i,w_i$ have the same form, the identifies~\eqref{id1}--\eqref{id4} and Lemma~\ref{Lem-m-1} carry through. Furthermore, $LR_\beta=My^kE_\beta M^{-1}$ and therefore $LR_\beta$ can be found from the right hand side of~\eqref{Lk} by substituting $y^k\to y^k E_\beta$. After that,repeating the same steps as in the proof of Lemma~\ref{L-}, we obtain that
\begin{multline}
 \big(My^kE_\beta M^{-1}\big)_-=\big(L^kR_\beta\big)_-=\sum_{i,j}\epsilon_{i}E_\beta w_{i-k}Y^k\wy^{-1}\wx^{-1}v_j\epsilon_j \\
-\sum_{i,j}\epsilon_iw_i\wx^{-1}\wy^{-1}Y^kv_{j+k}E_\beta\epsilon_{j}
-\sum_{i,j,\ell}\epsilon_iw_i\wx^{-1}\wy^{-1}Y^k v_\ell E_\beta w_{\ell-k}\wy^{-1}\wx^{-1}v_j\epsilon_j\\+\sum_{i,j,\ell}\sum_{a=0}^{k-1}\epsilon_iw_i\wx^{-1}Y^a v_\ell E_\beta w_{\ell-k}Y^{k-1-a}\wy^{-1}\wx^{-1}v_j\epsilon_j\,. \label{LkminusNEmat}
\end{multline}
Note that the right-hand side of this expression is obtained from~\eqref{LkminusNE} by replacing
\begin{align}
 &w_{i-k}\to E_\beta w_{i-k}\,, &
 &v_{j+k}\to v_{j+k} E_\beta\,,
 &v_\ell w_{\ell-k}\to v_\ell E_\beta w_{\ell-k}\,. \label{repl}
\end{align}
It is then easy to check that the proof of Prop.~\ref{prop4.3} can be repeated by making the replacements~\eqref{repl} throughout. \qed

Now varying $\beta$ and $k$ in~\eqref{tFlowXYmat}, \eqref{tFlowVWmat}, we have a family of commuting Hamiltonian flows on $\CM_{\alpha,\lambda}$, so by the above Proposition we obtain solutions of the hierarchy.

\begin{theorem} \label{ThFlowsMat}
The map $\CM_{\alpha,\lambda}\to\Big(\mathop{\mathrm{Mat}}(d,\mathbb C)\otimes\P\Big)^{d+1}$ which assigns to $(X,Y,v_i,w_i)$ the elements $L$, $R_\alpha$ given by~\eqref{LM}, \eqref{mmat} sends the Hamiltonian flows~\eqref{tFlowXYmat}, \eqref{tFlowVWmat} to the flows of the matrix hierarchy~\eqref{tkbetaL}, \eqref{tkbetaRalpha}. Therefore solutions of the hierarchy can be found by integrating the flows~\eqref{tFlowXYmat}, \eqref{tFlowVWmat} with an arbitrary initial point in $\CM_{\alpha,\lambda}$.
\end{theorem}

\subsection{Darboux coordinates and integrable dynamics}
\label{secDarCoorMat}

For the full hierarchy, the solutions are parameterised by the varieties \eqref{cms} for the quiver $Q$ as in \eqref{q}. As in Section~\ref{secFullHier}, let us consider the Calogero--Moser space $\CM_{\alpha,\lambda}$ for the special choice of the dimension vector $\alpha=n\delta$ where $\delta=(1, \dots , 1)$. Let $\CM_{n,\lambda}$ denote this space; the formula~\eqref{dimqmat} the gives us that $\dim \CM_{n,\lambda}=2mnd$. The construction of Darboux coordinates on $\CM_{n,\lambda}$ will be done entirely analogously to Section~\ref{secFullHier}.

Consider a point of $\CM_{n,\lambda}$ represented by a collection of linear maps $X_i,Y_j,v_{i,\alpha}$, $w_{i,\alpha}$, assuming the operator $X_{m-1}\dots X_1X_0\in\End(V_0)$ has distinct non-zero eigenvalues $\nu_1^m,\ldots,\nu_n^m$. Such points form a dense subset $\CM'_{n,\lambda}\subset \CM_{n,\lambda}$. As in the case $d=1$, there exists bases in each of $V_i$ diagonalizing the matrices $X_i$ as follows:
\begin{equation}
\label{slicegmat}
X_i=\mathrm{diag}({\nu}_1, \dots, {\nu}_n)\,.
\end{equation}
Substituting this into the relations~\eqref{XYvwkmatalpha} we obtain the following formulas:
\begin{align}
\label{slice2ggmat}
&(Y_i)_{jj}=\mu_j+\frac{1}{\nu_j}\sum_{r=0}^{i}\big(\lambda_r+(\psi_j\varphi_j)_{rr}\big)+\frac{1}{\nu_j}\sum_{r=0}^{m-1}\frac{r-m}m \big(\lambda_r+(\psi_j\varphi_j)_{rr}\big), \\
 \label{slice3ggmat}
&(Y_i)_{jk}=-\sum_{l=0}^{m-1}\,\frac{(\psi_k\varphi_j)_{i-l,i-l}\,\nu_j^{l}\nu_k^{m-l-1}}{\nu_j^m-\nu_k^m}\qquad\qquad (j\ne k), \\
&(v_{i,\alpha})_j=(\varphi_j)_{\alpha i}\,,
 \qquad\qquad (w_{i,\alpha})_j=(\psi_j)_{i\alpha}
 \qquad\qquad (\alpha=1,\ldots,d), \label{sliceVWmat}
\end{align}
where $i=0,\ldots,m-1$, and $\varphi_j$ and $\psi_j$ are matrices of size $d\times m$ and $m\times d$, respectively, satisfying the constraint $\tr(\varphi_j\psi_j) = -|\lambda|$ ($j=1,\ldots,n$). As above, with each pair $(\varphi_j, \psi_j)$ we can associate a point $q_j$ of the symplectic variety
\begin{equation*}
Q_{m,d}=\big\{(\varphi, \psi)\in\Hom(\C^m,\C^d)\oplus\Hom(\C^d,\C^m)\mid \tr(\varphi\psi)=-|\lambda|\big\}/\C^\times,
\end{equation*}
where the group $\C^\times$ acts by $(\varphi,\psi)\mapsto(\kappa\varphi,\kappa^{-1}\psi)$.
Thus we obtain a local parameterisation of $\CM_{n,\lambda}'$ by $\{\nu_j, \mu_j, q_j\}$.

Again, in the case $j\ne k$ we can rewrite the formula as
\begin{equation}\label{slice4ggmat}
(Y_i)_{jk}=\sum_{r\in I} \,\frac{\mu^{-ir}\tr\left(\varphi_j \, E^r \,\psi_k\,\right)}{m(\nu_k-\mu^{r}\nu_j)}\,,\qquad E^r=\mathrm{diag}(1, \mu^{r}, \mu^{2r}, \dots, \mu^{(m-1)r}).
\end{equation} 

The symplectic structure~\eqref{symp1mat} induces a symplectic structure on $\CM_{n,\lambda}$, which we denote by $\omega_{\CM_{n,\lambda}}$.

\begin{prop}
The coordinates $\nu_j, \mu_j, q_j=(\varphi_j, \psi_j)$ are local Darboux coordinates on $\CM'_{n,\lambda}$, that is the symplectic form $\omega_{\CM_{n,\lambda}}$ can be written as
\begin{equation} \label{omegaCnmat}
 \omega_{\CM_{n,\lambda}}=\sum_{j=1}^n \big(m\,\mathrm{d}\mu_j \wedge \mathrm{d}\nu_j+\tr(\mathrm{d}\psi_j \wedge \mathrm{d}\varphi_j)\big)\,,
\end{equation} 
where $\tr(\mathrm{d}\psi_j \wedge \mathrm{d}\varphi_j)=\sum\limits_{r=0}^{m-1}\sum\limits_{\beta=1}^d\mathrm{d}(\psi_j)_{r\beta} \wedge \mathrm{d}(\varphi_j)_{\beta r}$ is understood as the corresponding form on $Q_{m,d}$.
\end{prop}

The derivation of~\eqref{omegaCnmat} is completely analogous to the proof of Prop.~\ref{sfs}.

Let us consider the variety $\hreg \times \mathfrak h^*\times (Q_{m,d})^n$ equipped with the symplectic form~\eqref{omegaCnmat}. Substituting \eqref{slice2ggmat}--\eqref{sliceVWmat} into the Hamiltonians $H_{k,\beta}$ we get Poisson commuting functions on this symplectic variety. As in the case $d=1$  we obtain an integrable system in this way:

\begin{prop} \label{liouv}
Let $H_{k,\beta}=H_{k,\beta}(\nu, \mu, q)$ be the functions of the variables $(\nu, \mu, q)\in \hreg \times \mathfrak h^*\times (Q_{m,d})^n$, obtained by substituting the formulas~\eqref{slice2ggmat}--\eqref{sliceVWmat} into $H_{k,\beta}=-\sum_{i\in I}\tr(E_\beta w_i Y^k v_{i+k})$. Then $H_{k,\beta}$ with $k=1, \dots, mn$ and $\beta=1,\ldots,d$ are functionally independent, Poisson commuting functions and hence define a completely integrable system. The Hamiltonian flows defined by $H_{k,\beta}$ extend from $\CM'_{n,\lambda}=\hreg \times \mathfrak h^*\times (Q_{m,d})^n/(S_n\ltimes\Gamma^n)$ to produce complete flows on the Calogero--Moser space $\CM_{n,\lambda}=\CM_{n\delta, \lambda}$ for the quiver \eqref{qd}.  
\end{prop}

\noindent{\bf Proof.} The number of integrals is $mnd$. It is the half of the dimension of the phase space $\hreg \times \mathfrak h^*\times (Q_{m,d})^n$. Thus it suffices to show their independence.

Note that there is an isomorphism $\CM_{n,\lambda}\to \CM_{n,\lambda'}$, where $\lambda'=(-\lambda_{m-1},\ldots,-\lambda_0)$, defined by
\begin{align*}
 &X_i\to Y_{-i-1}\,, &&Y_i\to X_{-i-1}\,,
 && v_{i,\alpha}\to v_{-i,\alpha}\,, && w_{i,\alpha}\to -w_{-i,\alpha}\,.
\end{align*}
Under this isomorphism the Hamiltonians $H_{k,\beta}$ takes the form
$$\wt H_{k,\beta}=\sum_{i\in I}w_{i+k,\beta} X^k v_{i,\beta}=\sum_{i\in I}\sum_{j=1}^n (\varphi_j)_{\beta,i} (\psi_j)_{i+k,\beta} \nu_j^k \,.$$
Therefore, it is enough to proof the independence of the functions $\wt H_{k,\beta}$ in the neighbourhood of a point of $\hreg \times \mathfrak h^*\times (Q_{m,d})^n$.

Consider $\wt H_k:=\sum\limits_{\beta=1}^d\wt H_{k,\beta}$ for $k=mp$, $p=1,\ldots,n$:
\begin{align*}
 \wt H_{mp}=\sum_{j=1}^n \tr(\varphi_j\psi_j)\nu_j^{mp}=-|\lambda|\sum_{j=1}^n\nu_j^{mp}\,.
\end{align*}
Since $\nu_i^m\ne\nu_j^m\ne0$ for $i\ne j$ on $\hreg$, the variables $\nu_1,\ldots,\nu_n$ can be expressed locally via $\wt H_{mp}$. Now let $k=mp-q$, where $p=1,\ldots,n$ and $q=0,\ldots,m-1$. Then
\begin{align*}
 &\wt H_{mp-q,\beta}=\sum_{j=1}^n z_{j,\beta,q}\nu_j^{mp-q}\,, &&\text{where}\quad
 z_{j,\beta,q}=\sum_{i\in I}(\varphi_j)_{\beta,i}(\psi_j)_{i-q,\beta}\,.
\end{align*}
Fixing $\beta$ and $q$ we obtain a system of linear equations with a Vandermonde matrix. It follows that locally $z_{j,\beta,q}$ are functions of $\wt H_{k,\beta}$, $k\le mn$. Thus it suffices to show that there are $mnd$ independent functions on $\hreg \times \mathfrak h^*\times (Q_{m,d})^n$ among $\nu_j$ and $z_{j,\beta,q}$, $j=1,\ldots,n$, $\beta=1,\ldots,d$, $q=0,\ldots,m-1$.

For any fixed $j$, the functions $z_{j,\beta,0}$ satisfy the constraint $\sum_{\beta=1}^d z_{j,\beta,0}=-|\lambda|$. Clearly, we only need to prove that among $z_{j,\beta,q}$ with $\beta=1,\ldots,d$, $q=0,\ldots,m-1$ there exist $md-1$ functions which are independent as functions on $Q_{m,d}$.

Let us first consider $z_{j,\beta,q}$ as functions on a lager space $\wh Q_{m,d}=\Hom(\C^m,\C^d)\oplus\Hom(\C^d,\C^m)$. Choose the point $P\in\wh Q_{m,d}$ with the following coordinates:
\begin{align*}
 &(\varphi_j)_{\beta,i}=-|\lambda|/d\,, &&(\psi_j)_{i,\beta}=\delta_{i,0}\,.
\end{align*}
Since
 $\dfrac{\partial z_{j,\beta,q}}{\partial(\varphi_j)_{\alpha,r}}=\delta_{\alpha\beta}\delta_{rq}$ at the point $P$, the differentials $dz_{j,\beta,q}$ ($\beta=1,\ldots,d$, $q=0,\ldots,m-1$) are linearly independent elements of $T^*_P\wh Q_{m,d}$. Note that $P$ belongs to the hypersurface $H\subset\wh Q_{m,d}$ defined by the equation $\sum_{\beta=1}^d z_{j,\beta,0}=-|\lambda|$. Therefore, when restricted to $T_P H$, these differentials will span a subspace in $T^*_P H$ of dimension $(md-1)$. It is clear that the same will be true in a neighbourhood of $P$ in $H$. This means that among $z_{j,\beta,q}$ there are $md-1$ elements which are independent as functions on $H$. It follows that they are also independent as functions on $Q_{m,d}=H/\C^\times$. \qed

It is convenient to consider the Hamiltonians
\begin{align*}
 H_k:=\sum\limits_{\alpha=1}^d H_{k,\alpha}=-\sum_{i\in I}\tr(w_i Y^k v_{i+k})\,.
\end{align*}
They are polynomials of the order $k$ in the momentum variables $\mu_i$. For $k=mp$ their leading terms are $H_{mp}=|\lambda|\sum_{j=1}^n\mu_j^{mp}$.
Note that $H_{mp}$ are functions of $\mu_j,\nu_j$ and $(\psi_k\varphi_j)_{ii}$, where $i\in I$ and $j,k=1,\ldots,n$.

\begin{example} \normalfont
 In the case $m=1$ the first two Hamiltonians have the form
\begin{align} \label{H1H2}
 &H_1=|\lambda|\sum_{i=1}^n\mu_i\,, &
 &H_2=|\lambda|\sum_{i=1}^n\mu_i^2+|\lambda|\sum_{i\ne j}\frac{(\psi_i\varphi_j)(\psi_j\varphi_i)}{(\nu_i-\nu_j)^2}\,.
\end{align}
This is the spin Calogero--Moser system from~\cite{GH}. Its relation with the matrix KP hierarchy and the spaces $\CM_{n,\lambda}$ has been considered in~\cite{KBBT}, \cite{W2}, \cite{BGK}, \cite{BP}. Note that even in the $m=1$ case, our Theorem~\ref{ThFlowsMat} gives a little bit more compared to~\cite[Theorem~3.5]{BGK}, since we consider the full multicomponent KP hierarchy, while in~\cite{BGK} only the flows corresponding to $\partial_{t_k}=\sum_{\beta=1}^d\partial_{t_{k,\beta}}$ are considered.
\end{example}

\subsection{Spherical case}
\label{secSphCase}

For the multicomponent hierarchy~\eqref{tkbetaL}, \eqref{tkbetaRalpha} we can consider a sub-hierarchy consisting of the equations
\begin{align}
 &\frac{\partial L}{\partial t_{mj,\beta}}=\big[(L^{mj}R_\beta)_+,L\big]\,, &
 &\frac{\partial R_\alpha}{\partial t_{mj,\beta}}=\big[(L^{mj}R_\beta)_+,R_\alpha\big] \label{tmjbeta}
\end{align}
($j\in\mathbb N_0$, $\beta=1,\ldots,d$) where the elements $L$ and $R_\alpha$ are assumed to have the following symmetry
\begin{align} \label{equivmat}
 &\epsilon_i L= L\epsilon_{i-1}\,, &&\epsilon_i R_\alpha=R_\alpha \epsilon_i\,.
\end{align}
These are analogous to the conditions~\eqref{sym}, and we will call this the {\it multicomponent spherical sub-hierarchy}.

Once again, solutions of this sub-hierarchy will be constructed by using representations of quivers.

In Section~\ref{secConSol} we constructed the Calogero--Moser space $\CM_{\alpha,\lambda}$ associated with the quiver~\ref{qd}. Let us now consider the quiver~\eqref{qld} and write $\CM^{sph}_{\alpha,\lambda}$ for the corresponding Calogero--Moser space. For simplicity, we restrict to the case $\ell=0$, so we have $d$ arrows $b_{0,1},\ldots,b_{0,d}$, all going from $0$ to $\infty$. In this case the variety $\CM^{sph}_{\alpha,\lambda}$ can be viewed as a symplectic subvariety of $\CM_{\alpha,\lambda}$ by setting $v_i=w_i=0$ ($i\ne0$) in the formulas~\eqref{XYvwkmatalpha}--\eqref{XYvwkmat}. By Prop.~\ref{dim}
\begin{align*}
 \dim \CM^{sph}_{\alpha,\lambda}=2d\alpha_0-\sum_{i\in I}(\alpha_i-\alpha_{i+1})^2\,.
\end{align*}

Now for $(X,Y,v_0,w_0)\in\CM^{sph}_{\alpha,\lambda}$ consider the element
\begin{align}
 &M=1-\epsilon_0w_0(X-x)^{-1}(Y-y)^{-1}v_0\epsilon_0\,. \label{Mmateq}
\end{align}
Since it commutes with each of $\epsilon_i$, the operators
\begin{align}
L=MyM^{-1}\,,\qquad\qquad R_\alpha=M E_\alpha M^{-1} \label{LRmatsph}
\end{align}
will have the symmetry~\eqref{equivmat}.


\begin{theorem}
Let $X$, $Y$, $w_0$ and $v_0$ satisfy the equations~\eqref{tFlowXYmat}, \eqref{tFlowVWmat} with $k=mj$ (where $v_i=w_i=0$ for $i\ne0$). Then the elements $L$, $R_\alpha$ given by~\eqref{Mmateq}, \eqref{LRmatsph} satisfy the equations of the multicomponent spherical sub-hierarchy~\eqref{tmjbeta}. Therefore solutions of the hierarchy can be found by integrating the flows~\eqref{tFlowXYmat}, \eqref{tFlowVWmat} with an arbitrary initial point in $\CM^{sph}_{\alpha,\lambda}$.
\end{theorem}

\noindent{\bf Proof.} The flows~\eqref{tFlowXYmat}, \eqref{tFlowVWmat} preserve the subspace $\CM^{sph}_{\alpha,\lambda}\subset\CM_{\alpha,\lambda}$ defined by the conditions $v_i=0$, $w_i=0$ ($i\ne0$). Therefore the proof reduces to the application of Prop.~\ref{prop6.4} in the case $k=mj$. \qed

The restriction of the flows~\eqref{tFlowXYmat}, \eqref{tFlowVWmat} with $k=mj$ onto $\CM^{sph}_{\alpha,\lambda}$ is given by the Hamiltonians
\begin{align*}
 &h_{j,\beta}=H_{mj,\beta}\big|_{\CM^{sph}_{\alpha,\lambda}}=-w_{0,\beta} Y^{mj}v_{0,\beta} &&(j\in\mathbb N_0,\; \beta=1,\ldots,d).
\end{align*}
Consider the principal Hamiltonians $h_j:=\sum_{\beta=1}^d h_{j,\beta}$. Note that $h_j=\frac{|\lambda|}m \tr Y^{mj}=|\lambda|\tr\wh Y^j$, where $\wh Y=Y_0Y_1\cdots Y_{m-1}$.

Let us now consider the special case $\alpha_0=\alpha_1=\ldots=\alpha_{m-1}=n$, writing $\CM^{sph}_{n,\lambda}$ for $\CM^{sph}_{n\delta,\lambda}$. This is a $2nd$-dimensional symplectic variety with coordinates $\{\nu_j,\mu_j,q_j\}$, where $q_j=(\varphi_j,\psi_j)\in Q_{1,d}$ with
\begin{align*}
Q_{1,d}=\big\{(\varphi,\psi)\in(\C^d)\oplus(\C^d)^*\mid\tr(\varphi\psi)=\psi\varphi=-|\lambda|\big\}/\C^\times\,,\quad\dim Q_{1,d}=2d-2\,.
\end{align*}
A generic point of $\CM^{sph}_{n,\lambda}$ in these coordinates is described by the formulas
\begin{align}
 &(X_i)_{jk}=\nu_j\delta_{jk}\,, \label{Xkij0mat} \\
 &(Y_i)_{jj}=\mu_j+\frac{1}{\nu_j}\Big(\sum_{r=0}^{i}\lambda_r-\sum_{r=0}^{m-1}\frac{m-r}{m}\lambda_r\Big)\,, \label{Ykii0mat} \\
 &(Y_i)_{jk}=-\frac{\nu_j^i\nu_k^{m-i-1}}{\nu_j^m-\nu_k^m}\psi_k\varphi_j\,, &&&&(j\ne k), \\
 &(v_{0\alpha})_j=(\varphi_j)_\alpha\,, \qquad\qquad (w_{0\alpha})_j=(\psi_j)_\alpha\,, \label{v0w00mat}
\end{align}
where $i=0,\ldots,m-1$. A calculation similar to the proof of Prop.~\ref{sfs} shows that the symplectic form on $\CM^{sph}_{\alpha,\lambda}$ looks as $\omega_{sph}=\sum_{j=1}^n \big(m\mathrm{d}\mu_j\wedge \mathrm{d}\nu_j+\mathrm{d}\psi_j\wedge \mathrm{d}\varphi_j\big)$. By substituting the formulas~\eqref{Xkij0mat}--\eqref{v0w00mat} into the above $h_{j,\beta}$ we obtain Poisson commuting functions $h_{j,\beta}(\nu,\mu,q)$ on $\hreg\times\mathfrak h^*\times (Q_{1,d})^n$. By repeating the proof of Prop.~\ref{liouv} with $q=0$ one proves that $h_{j,\beta}$ with $j=1, \dots, n$ and $\beta=1,\ldots,d$ are functionally independent. Thereby we obtain the analogous result:

\begin{prop} \label{liouvSph}
 The functions $h_{j,\beta}(\nu,\mu,q)$ define a completely integrable system on $\hreg\times\mathfrak h^*\times (Q_{1,d})^n$. The Hamiltonian flows defined by $h_{j,\beta}$ extend from $\hreg \times \mathfrak h^*\times (Q_{1,d})^n/(S_n\ltimes\Gamma^n)$ to produce complete flows on $\CM^{sph}_{n,\lambda}$.  
\end{prop}

\begin{example} \label{ExCMBnmat} \normalfont
 Consider the case $m=2$. The Hamiltonians $h_{p,\beta}$ describe a systems of $n$ particles with coordinates $\nu_j$ momenta $\mu_j$ and spin variables $\varphi_j=(\varphi_{j1},\ldots,\varphi_{jd})^\top$ and $\psi_j=(\psi_{j1},\ldots,\psi_{jd})$ satisfying $\psi_j\phi_j=-|\lambda|$. The first principal Hamiltonian of the system has the form
\begin{align*}
 h_1=|\lambda|\sum_{i=1}^n\mu_i^2-\frac{|\lambda|\lambda_1^2}{4}\sum_{i=1}^n\frac{1}{\nu_i^2}-\frac{|\lambda|}{2}\sum_{i< j}\Big(\frac{1}{(\nu_i-\nu_j)^2}+\frac{1}{(\nu_i+\nu_j)^2}\Big) (\psi_i\varphi_j)(\psi_j\varphi_i)\,.
\end{align*}
This is a $B_n$ version of the Gibbons--Hermsen Hamiltonian~\eqref{H1H2}, cf.~\cite{LX}.
\end{example}

\section{Further directions and links}

Our results raise a number of questions; below we mention some of them and indicate further connections.

\subsection{} There are many aspects of the integrability for the usual KP hierarchy: $\tau$-function formalism, $W$-symmetries, (bi-)Hamiltonian approach, free-field realisation, and so on. Therefore it is a natural problem to find analogous notions for our hierarchy and to interpret the constructed solutions from that point of view.

\subsection{} For special values of $m$ the Calogero--Moser system for $G=\Z_m\wr S_n$ admits natural trigonometric and elliptic versions. In such cases, we expect these systems to be related to certain generalisations of the KP hierarchy. There is also a natural $q$-analogue of this story related to multiplicative preprojective algebras of Crawley-Boevey and Shaw~\cite{CBS}. These questions will be a subject of future work.

\subsection{} There is a link between the Calogero--Moser spaces for the cyclic quiver and the families of bispectral differential operators of Bessel type from~\cite{DG}, \cite{BHY}. The key to that link is the fact that $\epsilon_\ell y^m\epsilon_\ell$ can be interpreted as a generalised Bessel operator (the classical Bessel operator corresponds to $m=2$). Then one can use the operator~\eqref{ml} as a dressing, and construct bispectral differential operators in the form $\epsilon_\ell Mp(y^m)M^{-1}\epsilon_\ell$ for suitable polynomials $p$. The details will appear elsewhere.

Note that M.~Rothstein suggested in~\cite{Ro} a description of bispectral operators of Bessel type in terms of Wilson's Calogero--Moser spaces $\CM_n$. It is likely that his formulas describe some part of our Calogero--Moser space $\CM_{n,\lambda}$, but it is difficult to see this directly because the formulas in~\cite{Ro} are rather cumbersome. Also, the link with the Calogero--Moser problem for $\Z_m\wr S_n$ was not observed there.

\subsection{} Our results can be used to prove a conjecture of Etingof and Rains. Namely, in~\cite{ER} the authors study algebraically integrable operators with $\Z_m$-symmetry. For the case of rational coefficients they conjecture that the poles of such operators should be equilibrium points of the Calogero--Moser system for $\Z_m\wr S_n$ (see~\cite[Conjecture~4.10]{ER}). This is related to our results as follows. Let $\CM_{n,\lambda}$ be Calogero--Moser space for the quiver~\eqref{ql} with $\ell=0$ and let $(X,Y,v_0,w_0)\in\CM_{n,\lambda}$ be an equilibrium point for the Hamiltonian $H_m$. Then we have $\partial_{t_m}M=0$, where $M$ is~\eqref{ml}, implying $(L^m)_-=0$. Thus $D=\epsilon_0L^m\epsilon_0$ gives a purely differential operator. One has to restrict to the case $k_i\in\mathbb Z\tau$ to ensure that the operator $D$ is algebraically integrable. Let us assume that the matrix $X$ has distinct non-zero eigenvalues $\{\mu^r x_j\}$, $j=1,\ldots,n$, $r=0,\ldots,m-1$, and so $(x_1,\ldots,x_n)\in\mathbb C^n$ is an equilibrium point of the Calogero--Moser problem for $\Z_m\wr S_n$ by our assumption about $(X,Y,v_0,w_0)$. One can check that the resulting operator $D$ satisfies all the assumptions of~\cite[Conjecture~4.10]{ER}. This way we construct an algebraically integrable operator with $\mathbb Z_m$-symmetry starting from any equilibrium point of the Calogero-Moser problem. The full proof of the conjecture requires showing that the above construction produces all such operators. We will return to this problem elsewhere.

\subsection{} In the case $m=1$ the Calogero--Moser spaces $\CM_n$ are known to be related to the Bethe algebra of the quantum Gaudin model, see~\cite{MTV}. It would therefore be interesting to relate the Calogero--Moser spaces $\CM_{n,\lambda}$ with some generalisation of the Gaudin model. Note that a cyclotomic version of the Gaudin model has been proposed recently in~\cite{VY}. 

\subsection{} Finally, it would be interesting to understand which of our results can be generalised to other quivers (see in particular Remark~\ref{remGenQuiver} above).

\bibliographystyle{amsalpha}

\end{document}